\documentclass[lettersize,journal, ]{IEEEtran}
\usepackage{amsmath,amsfonts}
\usepackage{algorithmic}
\usepackage{array}
\usepackage{textcomp}
\usepackage{stfloats}
\usepackage{url}
\usepackage{verbatim}
\usepackage{graphicx}
\def\BibTeX{{\rm B\kern-.05em{\sc i\kern-.025em b}\kern-.08em
    T\kern-.1667em\lower.7ex\hbox{E}\kern-.125emX}}
\usepackage{balance}
\usepackage{verbatim}
\usepackage{amsfonts}
\usepackage{amssymb}
\usepackage{stfloats}
\usepackage{cite}
\usepackage{soul}
\usepackage{setspace}
\usepackage{graphicx}
\usepackage{psfrag}
\usepackage{amsmath}
\usepackage{array}
\usepackage{epstopdf}
\usepackage{authblk}
\usepackage{graphicx} 
\usepackage{amsthm} 
\usepackage{lipsum}
\usepackage{verbatim} 
\usepackage{authblk}
\usepackage{mathtools}
\usepackage{cuted}
\usepackage[lined,boxed,ruled]{algorithm2e}
\usepackage{booktabs}
\usepackage{subfigure}
\usepackage[font=small]{caption} % Add this line for smaller captions

\usepackage{ifthen}

\usepackage[usenames]{color}

\newcommand{\red}{\color{red}}

\newtheorem{Lemma}{Lemma}
\newtheorem{Theorem}{Theorem}

\newtheorem{Proposition}{Proposition}
\newtheorem{Remark}{Remark}

\DeclareMathOperator*{\argmax}{argmax}
\DeclareMathOperator*{\argmin}{argmin}
\DeclareMathOperator{\diag}{diag}

\newboolean{showaircomp}
\setboolean{showaircomp}{False} % 设置条件为true

\begin{document}

\title{Ultra-Low-Latency Edge Inference \\ for Distributed Sensing}
\author{{Zhanwei~Wang,~\IEEEmembership{Graduate Student Member,~IEEE},
Anders~E.~Kalør,~\IEEEmembership{Member,~IEEE},
    You~Zhou,~\IEEEmembership{Graduate Student Member,~IEEE},
Petar~Popovski,~\IEEEmembership{Fellow,~IEEE},
    and Kaibin~Huang,~\IEEEmembership{Fellow,~IEEE}}%

\thanks{The first manuscript was submitted on July 18, 2024—symbolically inspired by the ultra-fast Porsche 718—to reflect our pursuit of ultra-low-latency edge inference :)
}
\thanks{
The work of Z. Wang, Y. Zhou, and K. Huang was supported in part by the Research Grants Council of the Hong Kong Special Administrative Region, China under a fellowship award (HKU RFS2122-7S04), NSFC/RGC CRS (CRS\_HKU702/24), the Areas of Excellence scheme grant (AoE/E-601/22-R), Collaborative Research Fund (C1009-22G), and the Grants 17212423 \& 17304925, and in part by the Shenzhen-Hong Kong-Macau Technology Research Programme (Type C) (SGDX20230821091559018).  
} 
\thanks{The work of A.~E.~Kal{\o}r was supported in part by the Independent Research Fund Denmark (IRFD) under Grant 1056-00006B and in part by the Horizon Europe SNS project ``6G-GOALS'' (grant 101139232). The work of P.~Popovski was supported by the Villum Investigator Grant ``WATER'' from the Velux Foundation, Denmark and in part by the Horizon Europe SNS project ``6G-GOALS'' (grant 101139232).}%
  \thanks{Z.~Wang, Y.~Zhou, and K. Huang are with the Department of Electrical and Electronic Engineering, The University of Hong Kong, Hong Kong (email: \{zhanweiw, zhouyeee, huangkb\}@eee.hku.hk).

  A.~E.~Kal{\o}r is with the Department of Information and Computer Science, Keio University, Yokohama 223-8522, Japan, and with the Department of Electronic Systems, Aalborg University, 9220 Aalborg, Denmark (e-mail: aek@keio.jp).
  P.~Popovski is with the Department of Electronic Systems, Aalborg University, Denmark (email: petarp@es.aau.dk). \emph{Corresponding author: K.~Huang.}
  }

  }

\maketitle
% \begin{singlespace}

\begin{abstract}
There is a broad consensus that artificial intelligence (AI) will be a defining component of the sixth-generation (6G) networks. As a specific instance, AI-empowered sensing will gather and process environmental perception data at the network edge, giving rise to integrated sensing and edge AI (ISEA). Many applications, such as autonomous driving and industrial manufacturing, are latency-sensitive and require end-to-end (E2E) performance guarantees under stringent deadlines. However, the 5G-style ultra-reliable and low-latency communication (URLLC) techniques designed with communication reliability and agnostic to the data may fall short in achieving the optimal E2E performance of perceptive wireless systems.
In this work, we introduce an ultra-low-latency (ultra-LoLa) inference framework for perceptive networks that facilitates the analysis of the E2E sensing accuracy in distributed sensing by jointly considering communication reliability and inference accuracy.
By characterizing the tradeoff between packet length and the number of sensing observations, we derive an efficient optimization procedure that closely approximates the optimal tradeoff.
We validate the accuracy of the proposed method through experimental results, and show that the proposed ultra-Lola inference framework outperforms conventional reliability-oriented protocols with respect to sensing performance under a latency constraint.
\end{abstract}


% ---------------------------------------------------------------
%  Section : Introduction
% 
% ---------------------------------------------------------------
\section{Introduction}

Two of the main novelties expected in the upcoming \emph{sixth-generation} (6G) of mobile networks are edge \emph{artificial intelligence} (AI) and distributed sensing~\cite{HUAWEI-6G-2022,zhiyan_ISEA_survey}. The former will feature ubiquitous deployment of AI and learning algorithms at the network edge to automate mobile applications \cite{GX-CM-2020, 6ggoals, huiling_BS_2025}, whereas distributed sensing will facilitate large-scale sensor networks with, e.g., environment perception and localization capabilities~\cite{cui2025quansensing,CHEN2025,cui2025multiband}. 
The natural convergence of edge AI and distributed sensing, termed \emph{integrated sensing and edge AI} (ISEA), will provide a platform that promises to revolutionize numerous Internet-of-Things domains, such as auto-driving, health-care, smart cities, and industrial manufacturing \cite{Sherman-CoST-2022, liu2025integratedsensingedgeai}.
However, aligning ISEA with the ambitious 6G goals for achieving unprecedented \emph{end-to-end} (E2E) low-latency, reliability, and mobility requires new breakthroughs in ISEA technologies. One representative application scenario is the Level 4/5 autonomy in autonomous driving that requires real-time object detection within 30 milliseconds and with an accuracy close to 100\% \cite{Challenges-Level4/5}. 
As a step towards addressing this problem, we present the novel framework of \emph{ultra-low-latency} (ultra-LoLa) Inference for ISEA, which enables ultra-LoLa \emph{short-packet transmission} (SPT) of features extracted from remote sensing data. By carefully optimizing the tradeoff between communication reliability and sensing quality, the framework maximizes the E2E sensing accuracy given a stringent task completion deadline.

% Enabled by \emph{Short-Packet Transmission} (SPT), the framework provides an accurate characterization of the trade-off between E2E classification accuracy and communication reliability under strict latency constraints, and allows for optimization of the trade-off to enable reliable ultra-LoLa distributed sensing.

% While AI and sensing are critical enablers for ISEA, many applications also impose strict latency and reliability constraints. One representative scenario is the Level 4/5 autopilot that requires real-time object detection within 50 milliseconds and with an accuracy close to 100\% \cite{Challenges-Level4/5}.
% As a step towards supporting such applications, in this work we present the novel framework for designing and analyzing protocols for \emph{ultra-low-latency} (ultra-LoLa) Inference for ISEA. Enabled by \emph{Short-Packet Transmission} (SPT), the framework provides an accurate characterization of the trade-off between E2E classification accuracy and communication reliability under strict latency constraints, and allows for optimization of the trade-off to enable reliable ultra-LoLa distributed sensing.

The development of ultra-LoLa transmission techniques accelerated in 5G where \emph{ultra-reliable and low latency communications} (URLLC) was defined as one of 5G missions to support machine-to-machine applications \cite{Petar-IEEEProc-2016}.
Relevant techniques aim at realizing E2E latencies below one millisecond and packet-error rates of $10^{-5}$ or lower, allowing mission-critical applications (e.g., industrial automation) to be deployed in 5G networks.
A main challenge confronting URLLC is to balance the tradeoffs among data rate, latency, and reliability.
The short packets (or blocklengths) imposed by the strict deadline introduce a non-negligible transmission error probability, which depends on the operating regime~\cite{Petar-IEEEProc-2016}. When the channel diversity is limited and the transmitter only has statistical channel information, the reliability is dominated by deep fades and can be accurately characterized by outage expressions for infinite-length packets. On the other hand, when the transmitter has channel information or can leverage channel diversity, which is the focus in this paper, the reliability is mainly due to noise-induced decoding errors in the receiver.
%This is in contrast to the traditional paradigm where the effect of noise is suppressed across long-packet transmissions.
The resulting tradeoffs between packet error rate, data rate, and code length is a classic topic studied by information theorists \cite{Verdu-TIT_2010}. 
Guided by these results, one approach that 5G engineers have adopted to implement URLLC is using SPT \cite{popovski2019wireless}. 
This comes at the cost of low data rates, which limits the applications of URLLC to those with relatively low-rate transmissions, e.g., transmission of commands to remote robots or data uploading by primitive sensors (e.g., humidity, temperature, or pollution)~\cite{Nallanathan-TWC-2020,zhao2023joint}. 
Nevertheless, researchers have made diversified attempts to alleviate the limitation by designing customized SPT techniques including non-coherent transmission \cite{Liva-TCOM-2019}, optimal framework structures \cite{Petar-TCOM-2017},  power control \cite{Quek-TCOM-2018},   wireless power transfer \cite{Schmeink-JSAC-2018}, and multi-access schemes \cite{Nallanathan-TWC-2020}.
%\pp{PP: We also need to mention that, when the main source of error is the unknown fading and thus unknown SNR at the transmitter, then the short packet effects disappear and we can comfortably work with outage expressions for infinite-length packets.}

% In the 6G era, URLLC falls short on meeting the requirements of ISEA in two aspects. 
However, the existing SPT approaches targeting URLLC are insufficient to meet the requirements for ISEA.
 This is because ISEA applications are relatively data-intensive and frequently require the transmission of high-dimensional features (e.g., $512\times 7\times 7 $ features extracted from ModelNet dataset \cite{ModelNet-Ref} using VGG16 \cite{simonyan2015deep}).
 The stringent latency requirements for real-time inference render the transmission of high-dimensional feature vectors using long packets impractical. This creates a communication bottleneck, as the system must simultaneously meet the latency constraints with small-payload short packets while supporting the high-rate demands of feature transmission.
Moreover, the E2E performance of sensing tasks (e.g., inference accuracy) cannot be accurately
quantified by decoding error probability, which is more suitable for communication reliability assessment.
While existing URLLC techniques designed for 5G
focus a ultra high packet-level reliability, the target of ISEA is to achieve maximum E2E accuracy under a stringent deadline constraint.
Hence, two challenges faced by the deployment of ISEA in the 6G systems are:
\emph{(i)} overcoming the communication bottleneck with ultra-LoLa feature uploading from distributed sensors;
\emph{(ii)} designing new SPT techniques based on E2E sensing performance metrics.

The proposed ultra-LoLa inference framework is based on the popular ISEA architecture, known as the \emph{multi-view convolutional neural network} (MVCNN), provisioned with wireless connections between distributed sensors and an edge server \cite{ModelNet-Ref}.
The sensors employ lightweight pre-trained neural networks to extract features from their locally gathered data, which are then sent to the server for aggregation and inference through a deep neural network model.
Diverse strategies have been proposed to tackle the communication bottleneck and enhance the E2E sensing performance.
The majority of existing works consider a simplified point-to-point system with a single sensor/transmitter, for which MVCNN reduces to an architecture known as split inference that divides a global model into two parts for on-device feature extraction and on-server inference \cite{Chen-TWC-2019, ZJ-CoM-2020}.
Unlike split learning \cite{lin2023efficient, lin2024hierarchical}, which emphasizes optimizing DNN parameters during the training phase, split inference is primarily concerned with forward propagation using well-trained models.
In this context, the optimal model splitting point can be adjusted according to bandwidth and latency constraints \cite{Chen-TWC-2019}, and the tradeoff between computational cost and communication overhead \cite{ZJ-CoM-2020}.
The communication overhead in ISEA can be further reduced via feature pruning and quantization \cite{ZJ-CoM-2020}, and task-relevant sensor scheduling \cite{who2com}.
Researchers have also studied the latency-accuracy tradeoffs in other ISEA strategies, such as parallel processing and early exiting \cite{Zhiyan-JASC-2023}, and designed DNN-based feature transmission using \emph{joint source-channel coding} (JSCC) to additionally cope with channel noise \cite{Gunduz-WCL-2023}.
From the perspective of multi-access for ISEA, researchers have constructed a class of simultaneous-access schemes, called  \emph{over-the-air computing} (AirComp), to be a promising solution for realizing over-the-air feature aggregation when there are many sensors \cite{wen2023task_2, Zhiyan-AirPooling, AirBreathing}.
Recent work has also explored the E2E design integrating radar, sensing, communication, and AirComp \cite{GX-TWC-2023,wen2023task_1}.
However, while AirComp brings many advantages, its uncoded analog transmission can expose data to interference and noise, and encounters the issue of compatibility with existing digital systems.
Overall, the existing wireless techniques for ISEA are based on the assumption of long-packet transmission, which is not valid in the ultra-LoLa regime required for ultra-LoLa ISEA, which motivates the current work.

In this work, we propose an ultra-LoLa inference framework designed for SPT to achieve the optimal E2E performance of ISEA systems, i.e., maximizing the sensing accuracy referring to the accuracy of inference on uploaded data features.
We consider an ultra-LoLa application with a stringent deadline (e.g., 1ms for 6G) such that SPT is needed.
Underpinning the proposed framework is an important tradeoff, as revealed in this work, that exists between packet length and number of views.
Specifically, longer packets ensure more reliable communication but result in a lower sensing accuracy due to fewer sensor views that can be uploaded before the deadline. This reliability-view tradeoff necessitates the optimization of packet length to maximize the E2E sensing accuracy. In particular, focusing only on reliability, as in traditional URLLC, results in a sub-optimal sensing performance.

We consider two sensing scenarios.
The first is \emph{multi-snapshot sensing}, where a single sensor captures a sequence of observations, such as images of a moving object, extracts features from each observation using a local pre-trained model. The features are then fused locally before being uploaded using point-to-point SPT to the server for inference.
The other is \emph{multi-view sensing}, where multiple sensors simultaneously capture different observations of the same object and upload their features over a multi-access channel to the edge server, which aggregates the features and performs joint inference.
We study the problem under the assumption that the feature distribution follows a \emph{Gaussian mixture model} (GMM), as widely used in machine learning (see, e.g., \cite{Yang-CVPR-2014,figueroa2019semi}).
With \emph{channel state information} (CSI) at the transmitter, channel inversion is leveraged to ensure a fixed \emph{signal-to-noise ratio}  (SNR) for SPT. This is practical frame-based cellular system, as it offers data containers of deterministic length, which is important for scheduling in multi-user scenarios.
Based on finite-blocklength information theory, 
we quantify the E2E sensing accuracy by the probability of correct classification of sensing objects, taking into account both sensing noise and transmission errors.

The key contributions and findings of our work are summarized as follows:
\begin{itemize}
    \item \textbf{E2E performance analysis:} Tractable E2E performance analysis of ISEA is known to be challenging. Relevant results do not exist for ISEA with SPT, i.e., ultra-LoLa inference. To tackle the challenge, we adopt the popular GMM as feature distribution and derive a tight lower bound on classification accuracy as a function of the number of object classes and the minimum pairwise discriminant gain. The bound is a function of the feature size per sensor view and the number of views. A communication-computation integrated approach is then adopted to combine the preceding result with finite-blocklength information theory~\cite{Verdu-TIT_2010}. As a result, tractable expressions for E2E sensing accuracy are derived for both multi-snapshot and multi-view sensing under a task-deadline constraint, revealing the inherent reliability-view tradeoff discussed earlier and confirming the importance of packet-length optimization.
    \item \textbf{Packet-length optimization:}  
    To facilitate E2E optimization, accurate surrogate functions are derived for the sensing accuracy, under both considered sensing scenarios. For multi-snapshot sensing, the surrogate is a \emph{log-concave} function of the product of successful transmission probability and classification accuracy.
    For multi-view sensing, the derived surrogate function depends on the expected number of successfully transmitted views within the deadline.
    The unimodality of the surrogate functions allows the problem of accuracy maximization under a deadline constraint to be solved efficiently. The resulting optimal packet length reflects the influence of the  SNR, sensing complexity (i.e., number of classes), and additional network parameters.
     \item \textbf{Experiments:} The analytical results are validated in ISEA experiments using both synthetic (i.e., GMM) and real datasets (i.e., ModelNet \cite{ModelNet-Ref}). The proposed ultra-Lola inference designs are shown to achieve new-optimal performance and furthermore outperform traditional reliability-centered URLLC techniques.
\end{itemize}

The remainder of the paper is organized as follows. Sensing, communication, and inference models are defined in Section \ref{Sec:ModelandMetrics}. The E2E performance analysis for ultra-LoLa inference is presented in Section \ref{Sec:E2E_Acc_Analysis}, followed by a discussion on the optimal packet length in Section \ref{Sec:PL-optimization}. Experimental results are given in Section \ref{Experimental_res}, and the paper is concluded in Section \ref{Conclusion}.

% ---------------------------------------------------------------
%  Section : System Model
% 
% ---------------------------------------------------------------

\section{Models and Metrics}\label{Sec:ModelandMetrics}

The considered ultra-LoLa edge inference system consists of a single edge server performing remote object classification using $K$ distributed sensors. Time is divided into recurring \emph{inference rounds} initiated by the edge server. In each round, the edge server requests the uploading of observations from sensors, which are assumed to capture the same object under different conditions (e.g., different angles and lighting conditions in the case of images). 
The uploaded observations are used for joint inference to recognize the associated object.

\subsection{Data and Inference Models}
Each of the sensors produces an observation of a common object class $\ell$, and feeds it to a local pre-trained model to generate a corresponding feature vector, $\mathbf{x}_{k}\in\mathbb{R}^N$. 
The feature vectors are assumed to be randomly drawn from a joint conditional distribution $(\mathbf{x}_{1},\ldots,\mathbf{x}_{K})\sim p_{\mathbf{x}_{1:K}}(\mathbf{x}_{1},\ldots,\mathbf{x}_{K}|\ell)$, and the common object class $\ell$ is assumed to be drawn uniformly from a set of $L$ classes, so that the joint feature distribution can be written as
\begin{equation}
\label{eq:joint_dis}
  (\mathbf{x}_{1},\ldots,\mathbf{x}_{K})\sim \frac{1}{L}\sum_{\ell=1}^Lp_{\mathbf{x}_{1:K}}(\mathbf{x}_{1},\ldots,\mathbf{x}_{K}|\ell).
\end{equation}
The distribution is next elaborated for two classifier models.

\subsubsection{Gaussian Mixture Model for Linear Classification}
For tractability, 
the subsequent analysis
focuses on the case where the feature vectors are drawn from a GMM~\cite{figueroa2019semi, zw2025AIoutage,zeng2024knowledge}. Specifically, we assume that each feature vector $\mathbf{x}_{k}$ is drawn independently from a Gaussian distribution with mean $\boldsymbol{\mu}_{\ell}\in \mathbb{R}^{N}$ and covariance matrix $\mathbf{C}\in \mathbb{R}^{N\times N}$. Note that the mean (or class centroid) depends on the class while the covariance matrix is assumed to be the same for all classes. Without loss of generality, we assume that $\mathbf{C}=\diag(C_{1,1},C_{2,2},\ldots,C_{N,N})$ is diagonalized, e.g., using singular value decomposition. The joint distribution in \eqref{eq:joint_dis} becomes
\begin{equation}
\label{data_dis}
(\mathbf{x}_{1},\ldots,\mathbf{x}_{K})\sim \frac{1}{L} \sum_{\ell=1}^L\prod_{k=1}^K\mathcal{N}\left(\mathbf{x}_{k} \mid \boldsymbol{\mu}_{\ell}, \mathbf{C}\right),
\end{equation}
where $\mathcal{N}\left(\mathbf{x}_{k} \mid \boldsymbol{\mu}_{\ell}, \mathbf{C}\right)$ denotes the Gaussian \emph{probability density function} (PDF) with mean $\boldsymbol{\mu}_{\ell}$ and covariance matrix $\mathbf{C}$. Finally, to characterize the discernibility between any pair of classes $\ell$ and $\ell'$ of the GMM, we define the discriminant gain as the symmetric \emph{Kullback-Leibler} (KL) divergence~\cite{zw2025AIoutage}:
\begin{equation}\small
\begin{split}
    g_{\ell,{\ell'}} &= \sf{KL}(\mathcal{N}(\boldsymbol{\mu}_\ell,\mathbf{C})\,||\, \mathcal{N}(\boldsymbol{\mu}_{\ell'},\mathbf{C}))+\sf{KL}(\mathcal{N}(\boldsymbol{\mu}_{\ell'},\mathbf{C})\,||\,\mathcal{N}(\boldsymbol{\mu}_\ell,\mathbf{C}))\\
    & =(\boldsymbol{\mu}_{\ell}-\boldsymbol{\mu}_{\ell'})^{\sf{T}}\mathbf{C}^{-1}(\boldsymbol{\mu}_{\ell}-\boldsymbol{\mu}_{\ell'}).
\end{split}
\end{equation}

We consider a \emph{maximum likelihood} (ML) classifier for the distribution in \eqref{data_dis}, which relies on a classification boundary between a class pair being a hyperplane in the feature space. 
Due to the uniform prior on the object classes, the ML classifier is equivalent to a \emph{maximum a posteriori} (MAP) classifier, and the estimated label $\hat{\ell}$ is obtained as
\begin{equation}
  \label{Maha_min_classifier}
    \begin{split}
        \hat{\ell}& =\argmax_{\ell} \log \Pr(\overline{\mathbf{x}}|\ell)\\
        & = \argmin_{\ell} z_\ell(\overline{\mathbf{x}}),
    \end{split}  
\end{equation}
where $z_\ell(\overline{\mathbf{x}})=(\overline{\mathbf{x}}-\boldsymbol{\mu}_\ell)^{\sf{T}}\mathbf{C}^{-1}(\overline{\mathbf{x}}-\boldsymbol{\mu}_\ell)$ is the squared Mahalanobis distance between the input feature vector $\overline{\mathbf{x}}$ and the class centroid $\boldsymbol{\mu}_{\ell}$.

\subsubsection{General Model for CNN Classification}
We also consider a more realistic but analytically intractable scenario in experiments, where feature vectors are extracted from images using a CNN, as shown in Fig. \ref{Fig:Diag_CNN}.
The neural network model comprises multiple convolutional layers followed by multiple fully connected layers and a softmax output activation function that outputs a confidence score of each label. The layers are split into a sensor sub-model and a server sub-model, represented as functions $f_{\sf sen}(\cdot)$ and $f_{\sf ser}(\cdot)$, respectively.
The feature vector is constructed by passing an image ${\mathbf{M}_k}$ of the common object through a pre-trained CNN, i.e., $\mathbf{x}_k=f_{\sf sen}(\mathbf{M_k})$. The feature vectors are then aggregated to $\overline{\mathbf{x}}$, and the confidence score of the server-side classifier can be obtained by feeding the aggregated feature vector into the server sub-model, i.e., $ \{s_1,\dots,s_\ell,\dots,s_L\} =f_{\sf ser}\left(\overline{\mathbf{x}}\right)$.
With such processing, the MVCNN classifier outputs the inferred label that has the maximum confidence score, i.e., $\hat{\ell}=\argmax_{\ell} s_\ell$.

\subsection{Sensing Models}

% \begin{figure*}[t]
% 	\centering	\begin{minipage}[b]{1\textwidth}
% 		\centering	\includegraphics[width=1\textwidth]{Figures/Sensing Method.pdf}
% 	\end{minipage}
% 	% }
% 	\caption{Diagrams of sensing models and CNN classifier in Ultra-LoLa edge inference system.}
%  % \AK{I suggest to simplify this figure and use it in the introduction instead of in the system model. Fig. 2 essentially shows the same thing (but less fancy).}
% 	\label{fig:system-diagram-1}
% % \vspace{-3mm}
% \end{figure*}

\begin{figure*}[t!]
\centering
\subfigure[Multi-snapshot sensing]{
\includegraphics[width=0.6\columnwidth]{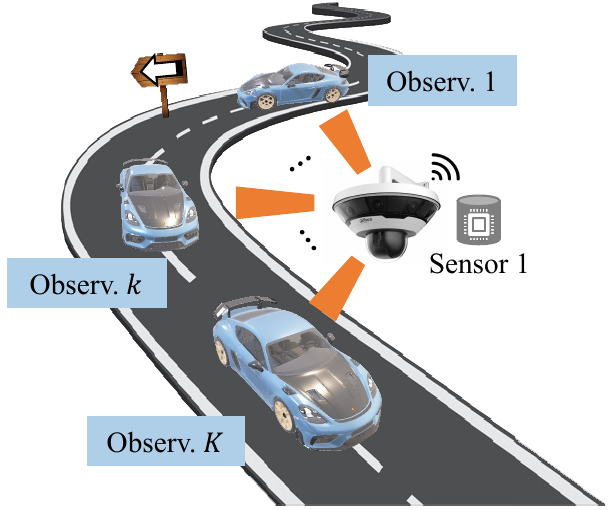}\label{Fig:Diag_MS}}
\subfigure[Multi-view sensing]{
\includegraphics[width=0.46\columnwidth]{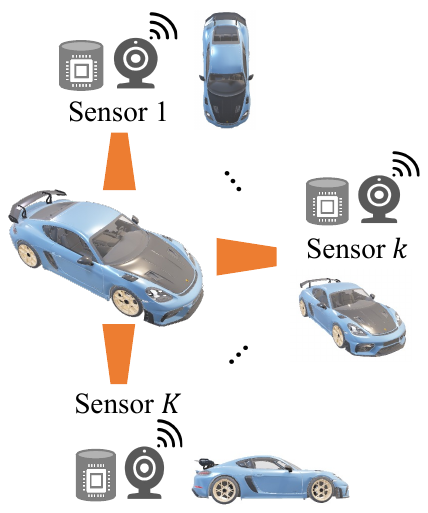}\label{Fig:Diag_MV}}
\subfigure[CNN classifier]{
\includegraphics[width=0.64\columnwidth]{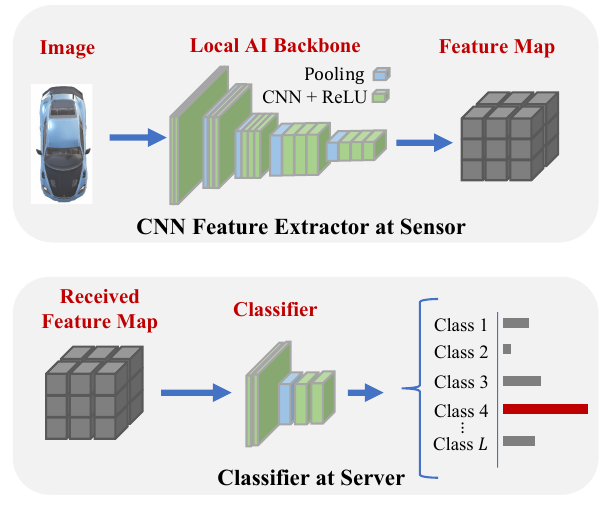}\label{Fig:Diag_CNN}}
\caption{Sensing models and CNN classifier in an ultra-LoLa edge inference system.}
\label{Fig:Diagram_1}\vspace{-3mm}
\end{figure*}

Two sensing scenarios as illustrated in Fig. \ref{Fig:Diagram_1} are modeled as follows.

\subsubsection{Multi-Snapshot Sensing}
As illustrated in Fig.~\ref{Fig:Diag_MS}, the $K$  observations are generated sequentially and then fused by a single device before being transmitted to the edge server. 
Multi-snapshot sensing represents, for instance, a scenario where the device captures a series of images of a moving vehicle in an intelligent traffic monitoring system. The temporal fusion of snapshots reduces motion blur and enhances resolution, enabling the sensor to transmit critical information to servers for inference \cite{biswas2016intelligent}.
We assume a sensing duration of $\Delta_S$ seconds, so that it takes a total of $T_{\sf{S}}=K\Delta_S$ seconds to generate the $K$ observations. The observations are then fused to an averaged feature vector $\overline{\mathbf{x}}\in\mathbb{R}^N$, given as
\begin{equation}
  \label{feature_avepool}
  \overline{\mathbf{x}}= \frac{1}{K}\sum_{k=1}^{K}\mathbf{x}_{k},
\end{equation}
which is transmitted to the edge server for inference.

\subsubsection{Multi-View Sensing}
Compared to multi-snapshot sensing by a single sensor, multi-view sensing is distributed where the $K$ observations are produced and transmitted in parallel by $K$ sensors, as shown in Fig.~\ref{Fig:Diag_MV}.
For instance, multi-view sensing provides comprehensive, multi-perspective data that reduces ambiguity, fills in occluded areas, and improves feature representation for more accurate and robust 3D shape reconstruction \cite{ModelNet-Ref}.
The total sensing duration is then simply $T_{\sf{S}}=\Delta_S$, and the edge server performs both feature fusion of the received feature vectors and inference.\footnote{We acknowledge the possibility of using AirComp to fuse views, but its lack of coding may not meet the high-reliability requirement that we target.}

\subsection{Short-Packet Transmission Model}
\label{Comm_model}
The feature vector at each sensor device is uploaded to the edge server over a wireless channel with a total bandwidth of $B_W$ Hz. Each feature vector is quantized to a sufficiently large number of bits $NQ_B$, such that the quantization error is negligible. The quantized features are then encoded using a channel code and transmitted over $D$ complex channel uses\footnote{As commonly done in practical systems, we assume that the payload of the packet (i.e., the features) are encoded separately from metadata. 
To simplify exposition and analysis, we neglect the impact of metadata (e.g., packet headers) on latency and error probability, noting that joint optimization of metadata and payload is a research topic on its own~\cite{popovski18urllc}. This assumption can be further justified by the fact that no channel-state information is required at the server, and thus scheduling can be pre-determined, e.g., by assigning periodical transmission slots.
}
(i.e., $D$ is the packet length and $NQ_B/D$ is the rate).

We assume random i.i.d.\ block-fading channels between the devices and the server, where the channel coefficient of device $k$ is $h_k$. The coefficient remains constant throughout the transmission of each feature vector. 
We assume that the devices have perfect knowledge of their channel coefficients via channel training, so that they can perform truncated channel inversion power control that is optimal in the finite blocklength regime~\cite{yang15powercontrol}. 
On the other hand, the server learns about the channels only after receiving the feature vector transmissions in the uplink, and thus cannot perform scheduling based on channel information. Following the principle of truncated channel inversion, device $k$ inverts its channel only if its channel gain $|h_k|^2$ exceeds a threshold $G_{th}$, and otherwise stays silent by setting its transmission power to zero.
\ifthenelse{\boolean{showaircomp}}{Unless otherwise stated, we}{We}
assume orthogonal transmissions such that the signal received by the server from device $k$ can be written as
\begin{equation}
\label{channel_model}
  \mathbf{y}_k = p_kh_k\mathbf{s}_k+\mathbf{w}_k,
\end{equation}
where $\mathbf{s}_k\in\mathbb{C}^D$ is the transmitted signal satisfying $(1/D)\|\mathbf{s}_k\|_2^2\le 1$, $\mathbf{W}_k\sim\mathcal{CN}(0,B_W N_0)$ is additive white Gaussian noise (AWGN) with noise spectral density $N_0$, and $p_k$ is the precoding coefficient given as
\begin{equation}
  \label{CI_PC}
  p_k=
  \begin{cases}
    \frac{\sqrt{P_{0}}h_k^*}{|h_k|^2},&|h_k|^2\geq G_{th}, \\
    0,&|h_k|^2< G_{th},
  \end{cases}
\end{equation}
where $P_{0}$ is the target signal power at the receiver.
The activation probability, denoted by $\xi_a$, is then obtained as 
\begin{equation}
  \label{active-prob}
  \xi_a=\Pr(|h_k|^2\geq G_{th}).
\end{equation}
We denote the subset of transmitting devices by $\mathcal{K}=\{k\in\{1,2,\ldots,K\} \mid |h_k|^2\ge G_{th}\}$.
The threshold $G_{th}$ is chosen such that $\mathbb{E}[|p_k|^2]=  \xi_a\mathbb{E}[|h_k|^{-2}] \le P_{\max}$, where $P_{\max}$ is the long-term power constraint of each sensor.

Under the above assumptions, whenever the channel gain exceeds $G_{th}$, the channel is effectively an AWGN channel with SNR $\gamma=\frac{P_0}{N_0 B_W}$~\cite{yang15powercontrol}. 
In this regime, the decoding error probability for a transmitting device ($k\in\mathcal{K}$) can be closely approximated as~\cite{Verdu-TIT_2010}
\begin{equation}
\varepsilon = Q\left(\ln (2) \sqrt{\frac{D}{V}}\left(\log _2(1+\gamma)-\frac{N Q_B}{D}\right)\right),\label{eq:varepsilon}
\end{equation}
where $Q(x)=\frac{1}{\sqrt{2 \pi}} \int_x^{\infty} e^{-\frac{t^2}{2}} d t$ and $V=1-(1+\gamma)^{-2}$ is the channel dispersion.
Combining this with the probability of outage (when $|h_k|^2<G_{th}$), the probability of successful transmission, denoted as $\rho$, can be expressed as
\begin{equation}
\label{succ_prob}
  \rho=\xi_a (1-\varepsilon).
\end{equation}
We define $\hat{\mathcal{K}}\subseteq \mathcal{K}$ as the set of sensors that successfully transmit their feature vector\footnote{We assume that packet errors can be detected perfectly by the server. In practice, this comes at a cost of a small reduction in the rate, e.g., caused by the inclusion of a cyclic redundancy check (CRC). However, we ignore this aspect to keep the presentation simple, noting that the added overhead (e.g., 8-24 bits in case of a CRC) would have only a small impact on the error probability for the blocklengths considered in this paper.}.

\subsection{System Configurations and Multi-Access}\label{sec:access_protocols}

%The system configurations and multi-access schemes for the two sensing scenarios are as follows.

\begin{figure*}[t!]
\centering
\subfigure[Multi-snapshot sensing with point-to-point transmission]{
\hspace{1cm}\includegraphics[width=0.65\columnwidth]{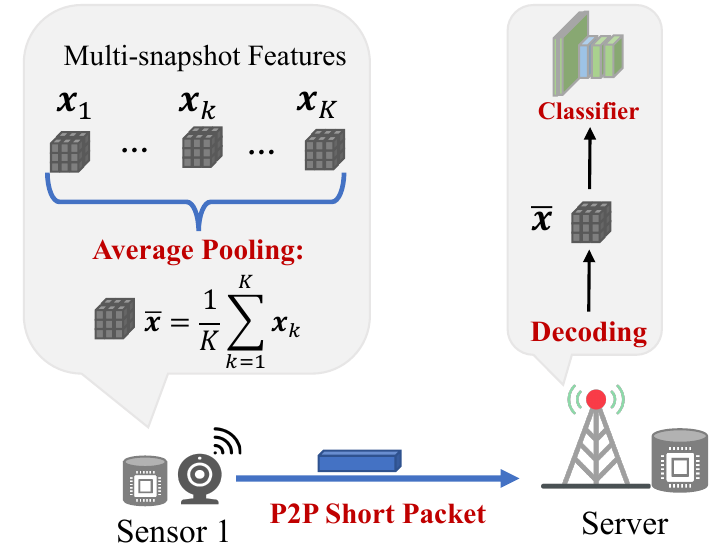}\label{Fig:Multi-access_P2P}\hspace{1cm}}
\subfigure[Multi-view sensing with orthogonal access]{
\hspace{1cm}\includegraphics[width=0.8\columnwidth]{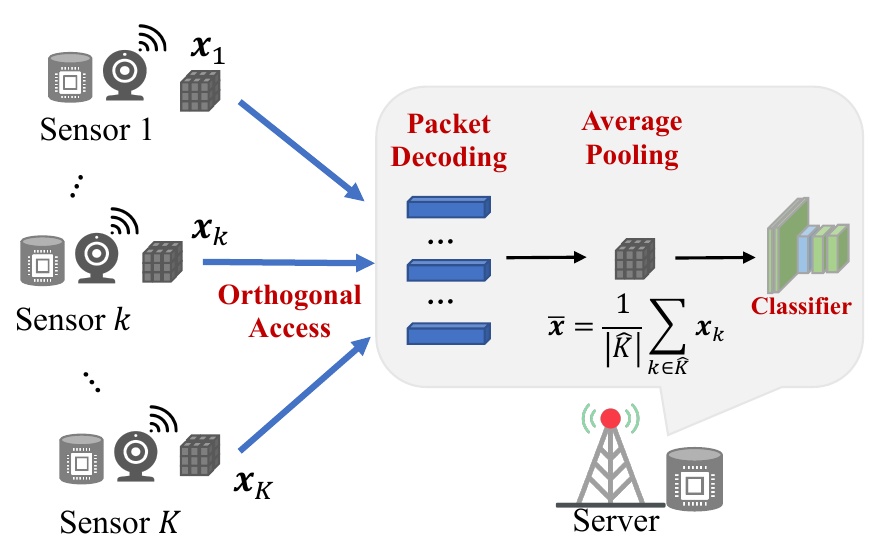}\label{Fig:Multi-access_OA}\hspace{1cm}}
\caption{Diagrams of multi-access models of the considered sensing tasks.}
\label{Fig:Diagram_2}\vspace{-3mm}
\end{figure*}

\subsubsection{Multi-Snapshot Sensing with Point-to-Point Transmission}
For multi-snapshot sensing, the quantized fused feature vector $\overline{\mathbf{x}}$ in \eqref{feature_avepool} is uploaded using the full system bandwidth $B_W$, as shown in Fig. \ref{Fig:Multi-access_P2P}. 
The transmission latency is then $T_{\sf{T}}=\frac{D}{B_W}$, where $D$ is the number of complex channel symbols.
In practical 5G protocols, the transmission latency can be controlled by configuring the number of symbols in each \emph{physical resource block} (PRB) \cite{3GPP_38_211}.
Upon successful transmission, the fused feature vector is used for inference at the edge server.

\subsubsection{Multi-View Sensing with Orthogonal Access}
In multi-view sensing, we consider \emph{time division multiple access} (TDMA) for feature uploading, as shown in Fig. \ref{Fig:Multi-access_OA}.
% Specifically, the system bandwidth is divided equally among the $K$ devices, so that a bandwidth $B=B_W/K$ is allocated to each device (regardless of their channels).
Specifically, the transmission of the packets is divided equally  into $K$ time slots, each of which has access to the full bandwidth of $B_W$.
Given the orthogonal transmission, each packet is decoded independently by the server with error probability $\varepsilon$ defined in \eqref{eq:varepsilon}. However, compared to the multi-snapshot setting the resulting transmission latency is $T_{\sf{T}}=\frac{DK}{B_W}$.
After the transmissions, the server performs average pooling over all successfully received feature vectors, i.e., it computes the averaged feature vector
\begin{equation}
\label{multiview_ap}\overline{\mathbf{x}}=\frac{1}{\hat{|\mathcal{K}|}}\sum_{k\in\hat{\mathcal{K}}}\mathbf{x}_{k},
\end{equation}
which is then fed into the classifier for inference.

\ifthenelse{\boolean{showaircomp}}{ 
\subsubsection{Multi-view Sensing with Over-the-Air Computing}
\label{MVS_Air}
Different from orthogonal access, \emph{over-the-air Computing} (AirComp) exploits the waveform superposition property \cite{RN30} to realize simultaneous multi-access.
Analog modulation for local features is required such that the received feature at server is the aggregated feature vector of scheduled devices.   
Considering the truncated channel inversion, the received signal is given as
\begin{equation}
\mathbf{y}= \sqrt{P_{0}}\sum_{k\in\mathcal{K}}\mathbf{x}_{k}+\mathbf{z},
\end{equation}
where $\mathbf{z}\sim \mathcal{CN}(\mathbf{0},N_0B_W \mathbf{I}_N) $ is the complex Gaussian vector representing channel noise from the analog transmission. Note that truncated channel inversion is not necessarily optimal for AirComp, but we apply it for analytical tractability and to allow for direct comparison between AirComp and digital transmission.
% Considering i.i.d channel coefficients over devices, the number of activation devices subjects to binomial distribution, i.e., $|\Tilde{\mathcal{K}}|\sim  \mathcal{B}(K,\xi_a)$.
Then taking the average pooling over active devices, the input of the classifier is given as
\begin{equation}
\label{AirComp features}
\overline{\mathbf{x}}=\Re\left({\frac{\mathbf{y}}{|\mathcal{K}|\sqrt{P_0}}}\right)=\frac{1}{|\mathcal{K}|}\sum_{k\in{\mathcal{K}}}\mathbf{x}_{k}+\Tilde{\mathbf{z}}
 \end{equation}
where $\Tilde{\mathbf{z}}\sim \mathcal{N}(\mathbf{0},\frac{1}{2|\mathcal{K}|^2\gamma'} \mathbf{I})$ is the scaled channel noise from the average pooling; $\gamma'=\frac{P_0}{N_0B_W}$ is denoted as the receive SNR of feature uploading from each device \cite{AirBreathing}.
The resultant transmission latency of AirComp depends on the feature size, given as $T_{\sf{T}}=NB_W$.
In the end, $\overline{\mathbf{x}}$ is fed into the classifier for inference.
}

\vspace{-3mm}
\subsection{E2E Performance Metric}

We measure the E2E performance by sensing accuracy.
The metric is defined as the probability of correct classification by taking into account both the random observations and the transmission process. Specifically, for a given number of snapshots $K$, the sensing accuracy in the multi-snapshot setting, denoted as $ \overline{E}_{\sf{ms}}$, is defined as
\begin{equation}
\label{eq:Expected_Acc_Single_Sensor}
    \overline{E}_{\sf{ms}}= \frac{\rho}{L}\sum_{\ell=1}^L \Pr(\hat{\ell}=\ell|\ell,K) + \frac{1-\rho}{L}, 
\end{equation}
where $\Pr  (\hat{\ell}=\ell|\ell,K)$ refers to the probability of correct classification with $K$ observations given the ground-truth label $\ell$. 
Note that the  term $\frac{1-\rho}{L}$ in \eqref{eq:Expected_Acc_Single_Sensor}  represents the probability of decoding failure followed by randomly guessing the correct label in an $L$-class classification task.
Similarly, in multi-view sensing, the sensing accuracy, denoted as $\overline{E}_{\sf{mv}}$, is given by
\begin{equation}
\begin{split}
\label{Expected Accuracy}
\overline{E}_{\sf{mv}}=\frac{1}{L}\sum_{\ell=1}^L \mathbb{E}_{|\hat{\mathcal{K}}|}\left[\Pr(\hat{\ell}=\ell|\ell,|\hat{\mathcal{K}}|)\mid \ell,K\right],
\end{split}
\end{equation}
where the expectation is over the number of received feature vectors $|\hat{\mathcal{K}}|$, taking into account that some of the $K$ users may fail to transmit their feature vectors.

\section{E2E Accuracy Analysis for Ultra-LoLa Inference
}
\label{Sec:E2E_Acc_Analysis}
In this section, we analyze the sensing accuracy of ultra-LoLa edge inference for the two considered sensing scenarios. 
% We ensure that the stringent task deadline is met by imposing an E2E latency constraint $T$ on the task completion of capturing images and feature transmission.
We consider an E2E deadline so that the sensing and transmission must be completed within a duration of $T$ seconds  (in the order of a few milliseconds). Specifically, given sensing time $T_{\sf{S}}$ and  transmission time $T_{\sf{T}}$, we require $T_{\sf{S}}+T_{\sf{T}}\le T$.

\subsection{Sensing Accuracy for Multi-Snapshot Sensing}

The computation of sensing accuracy in \eqref{eq:Expected_Acc_Single_Sensor} relies on the explicit expression of classification accuracy with a given number of observations $K$, i.e., 
$\frac{1}{L}\sum_{\ell=1}^L\Pr(\hat{\ell}=\ell|\ell,K)$.
However, it is difficult to obtain its closed-form expression due to the coupling among pairwise decision boundaries in a multi-class classifier. Instead, for a tractable analysis, we study the classification accuracy through the bound in Lemma \ref{inference_accuracy_LB}.

\begin{Lemma}[Lower Bound on Classification Accuracy]
\label{inference_accuracy_LB}
    Conditioned on successful transmission and number of classes $L\geq 2$,  the classification accuracy with $K$ snapshots, denoted as $A$, is lower bounded as 
\begin{align}
\small
     A&=\frac{1}{L}\sum_{\ell=1}^L \Pr(\hat{\ell}=\ell|\ell,K)\\
     &\geq \left[1-(L-1)Q\left( \frac{\sqrt{K g_{\min}}}{2}\right)\right]_{(1/L)^+},\label{eq:inference_accuracy_LB}
\end{align} 
where $[x]_{(1/L)^+}=\max\{1/L,x\}$ ensures an accuracy larger than random guessing probability $\frac{1}{L}$ and $g_{\min}=\min\{ g_{\ell,\ell'}|\ell,\ell'\in \{1,2,\dots,L\}, \ell\neq\ell' \}$ is the minimum pairwise discrimination gain among all classes.
\end{Lemma}
\noindent The proof is given in Appendix \ref{Closed_form_IA}.

\begin{Remark}[Classification Accuracy of Binary Classifier]
\label{Rem:Binary_Accuracy}
For a binary classifier with $L=2$, the classification accuracy in \eqref{eq:inference_accuracy_LB} reduces to the exact expression
    \begin{equation}
    \label{eq:binary_CCP}
        A=Q\left( -\frac{\sqrt{K g_{1,2}}}{2}\right),
    \end{equation}
    where $g_{1,2}=\sum_{n=1}^N\frac{(\boldsymbol{\mu}_1(n)-\boldsymbol{\mu}_2(n))^2}{C_{n,n}}$ is the discrimination gain.
\end{Remark}

As observed in Lemma \ref{inference_accuracy_LB} and also \eqref{eq:binary_CCP}, the lower bound on classification accuracy increases with the minimum discrimination gain and decreases with an increasing number of classes.
These two factors determine the difficulty of the classification task.

For a given packet length $D$, the maximum accuracy is obtained when the number of observations $K$ is maximized while satisfying the task deadline $T$. Considering the constraint $T_{\sf{S}}+T_{\sf{T}}\le T$ where $T_{\sf{S}}=K\Delta_S$, this quantity is given as
\begin{equation}
\label{eq:num_observations}
K_{\sf{ms}}=\left\lfloor\frac{T}{\Delta_S}-\frac{D}{B_W \Delta_S}\right\rfloor.
\end{equation}
Substituting \eqref{eq:inference_accuracy_LB} and \eqref{eq:num_observations} into \eqref{eq:Expected_Acc_Single_Sensor},  the  achievable sensing accuracy for multi-snapshot sensing is lower bounded in Proposition \ref{Proposition:E2E_Accuracy_SS}.

\begin{Proposition}[Sensing Accuracy for Multi-snapshot Sensing]
\label{Proposition:E2E_Accuracy_SS} 
For a given packet length $D$, the achievable sensing accuracy for multi-snapshot sensing with task completion deadline $T$ can be lower bounded as
\begin{equation}
\label{eq:E2E_Acc_SingleSensor}
\begin{split}
   \overline{E}_{\sf{ms}} 
   &\geq (L-1)\rho\left[ \frac{1}{L}- Q\left(\frac{\sqrt{K_{\sf{ms}}g_{\min}}}{2}\right)  \right]_{(1/L)^+}+\frac{1}{L},
   \end{split}
\end{equation}
where $[x]_{(1/L)^+}=\max\{1/L,x\}$, $\rho$ is the successful transmission probability given in \eqref{succ_prob}, and $K_{\sf{ms}}$ is given in \eqref{eq:num_observations}.
 \end{Proposition}

\begin{Remark}[Reliability-View Tradeoff]
    The lower bound on sensing accuracy in Proposition \ref{Proposition:E2E_Accuracy_SS} highlights the tradeoff between communication reliability and the number of views.  In particular, increasing the packet length $D$ increases the successful decoding probability, which in turn results in a higher successful transmission probability $\rho$, and vice versa.
    However, a large packet length comes at the cost of a smaller number of observations during the sensing phase, which decreases the discrimination gain $K_{\sf{ms}}G_{\min}$, and reduces the sensing accuracy due to increased sensing noise. 
Considering this tradeoff, the packet length $D$ should be optimized to maximize the sensing accuracy, which we do in Section \ref{Sec:PL-optimization}. 
\end{Remark}

\subsection{Sensing Accuracy for Multi-View Sensing}
Next, we investigate the reliability-view tradeoff for multi-view sensing by applying the lower bound on classification accuracy from Lemma \ref{inference_accuracy_LB} to the multi-view sensing accuracy in \eqref{Expected Accuracy}. Since each view is transmitted independently, the sensing accuracy for multi-view sensing, denoted as $\overline{E}_{\sf{mv}}$, can be lower bounded as
 \begin{equation}
\label{sdp+ms}
\begin{split}
% \overline{E}_{\sf{mv}} & \geq \mathbb{E}_{|\hat{\mathcal{K}}|}\left[ \left[1-(L-1)Q\left( \frac{\sqrt{|\hat{\mathcal{K}}| g_{\min}}}{2}\right)\right]_{(1/L)^+}  \right]\\
\footnotesize
\overline{E}_{\sf{mv}} & \geq  \hat{E}_{\sf{mv}}\\
&=\sum_{|\hat{\mathcal{K}}|=0}^{K} 
P_{|\hat{\mathcal{K}|}}
\left[1-(L-1)Q\left( \frac{\sqrt{|\hat{\mathcal{K}}| g_{\min}}}{2}\right)\right]_{(1/L)^+},
\end{split}
\end{equation}
where $\hat{\mathcal{K}}\subseteq \{1,2,\dots,K\} $ is the set of successfully transmitted packets, and $P_{|\hat{\mathcal{K}}|}= \binom{K}{|\hat{\mathcal{K}}|} \rho ^{|\hat{\mathcal{K}}|}(1-\rho)^{K-|\hat{\mathcal{K}}|}$  is the probability mass function (PMF) of the binomial distribution $|\hat{{\mathcal{K}}}|\sim\mathcal{B}(K,\rho)$.
Due to the TDMA transmissions, for a given task completion constraint $T\ge T_{\sf{S}}+T_{\sf{T}}$ and sensing duration $T_{\sf{S}}=\Delta_S$, the maximum number of sensors, denoted as $K_{\sf{mv}}$, is inversely proportional to $D$ and given as
    \begin{equation}
    \label{MS_K}
K_{\sf{mv}}=\left\lfloor\frac{(T-\Delta_S)B_W}{D} \right\rfloor.
    \end{equation}

\begin{figure}[t!]
\centering
\subfigure[Binary classification with $\gamma=5\text{dB}$]{
\includegraphics[width=0.45\columnwidth]{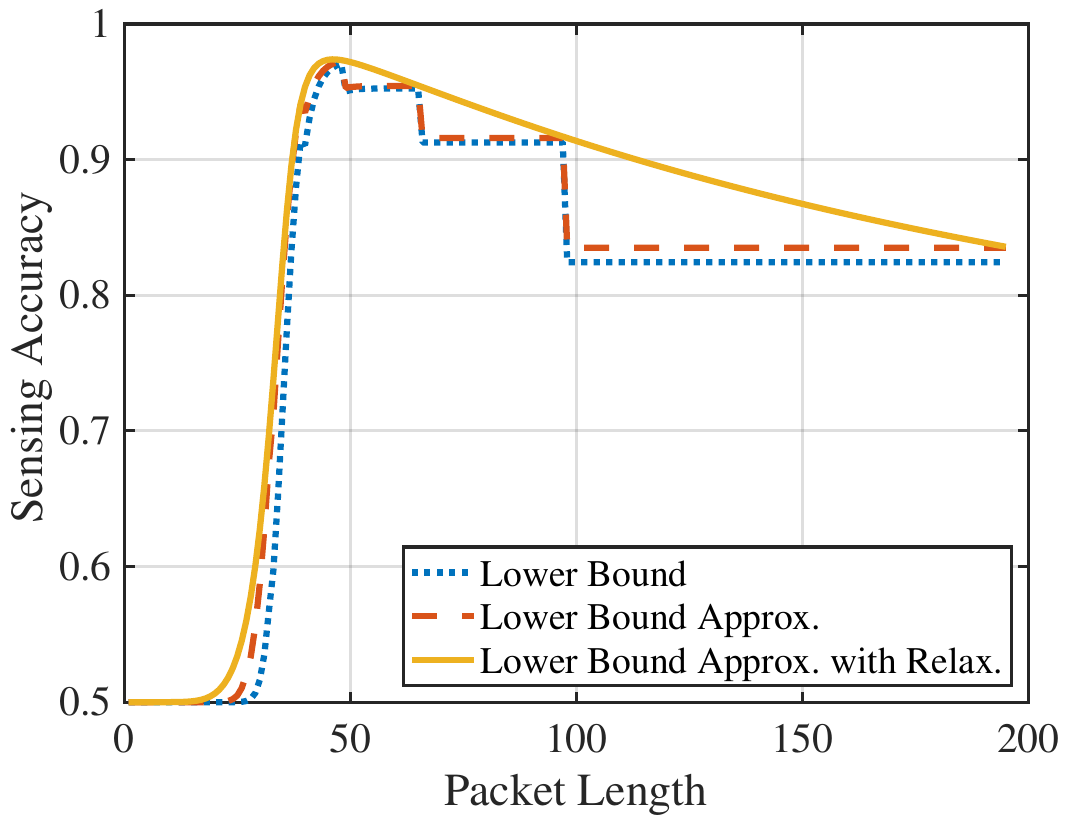}
\label{fig:NumRes_MV_PL_ACC(5dBL2)}
}
\subfigure[Binary classification with $\gamma=10\text{dB}$ ]{
\includegraphics[width=0.45\columnwidth]{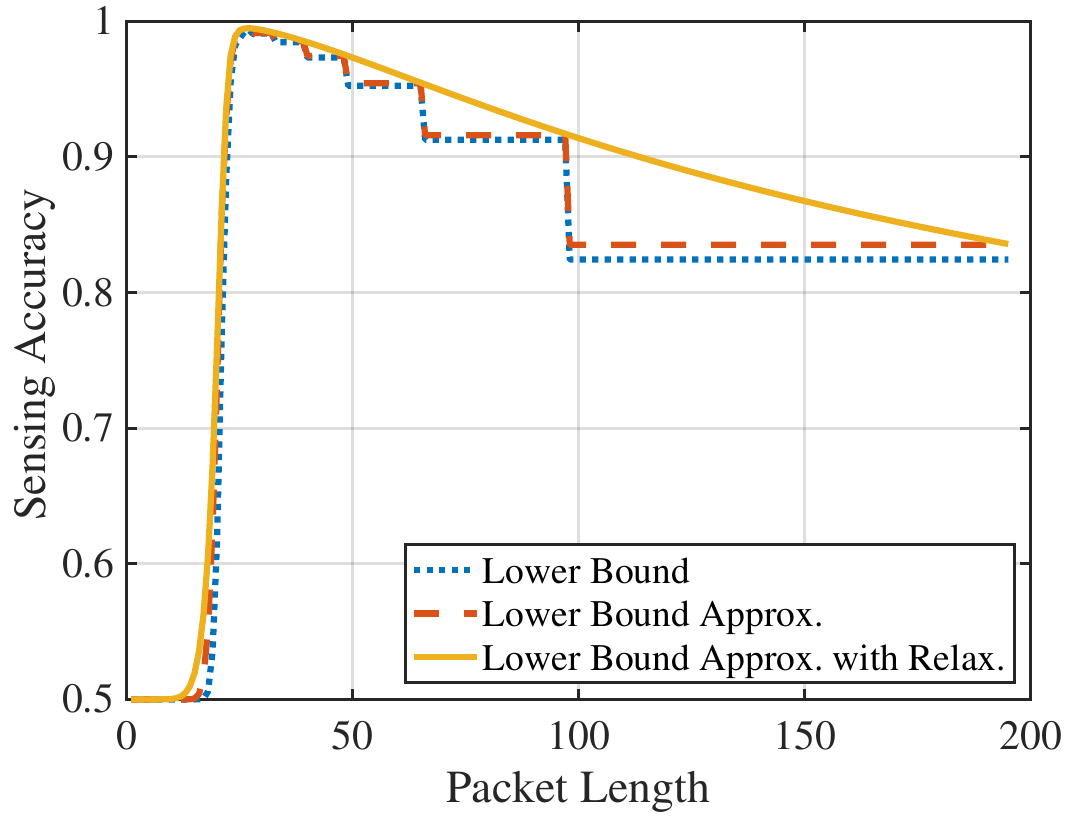}
\label{fig:NumRes_MV_PL_ACC(10dBL2)}
}
\subfigure[3 classes with $\gamma=5\text{dB}$ ]{
\includegraphics[width=0.45\columnwidth]{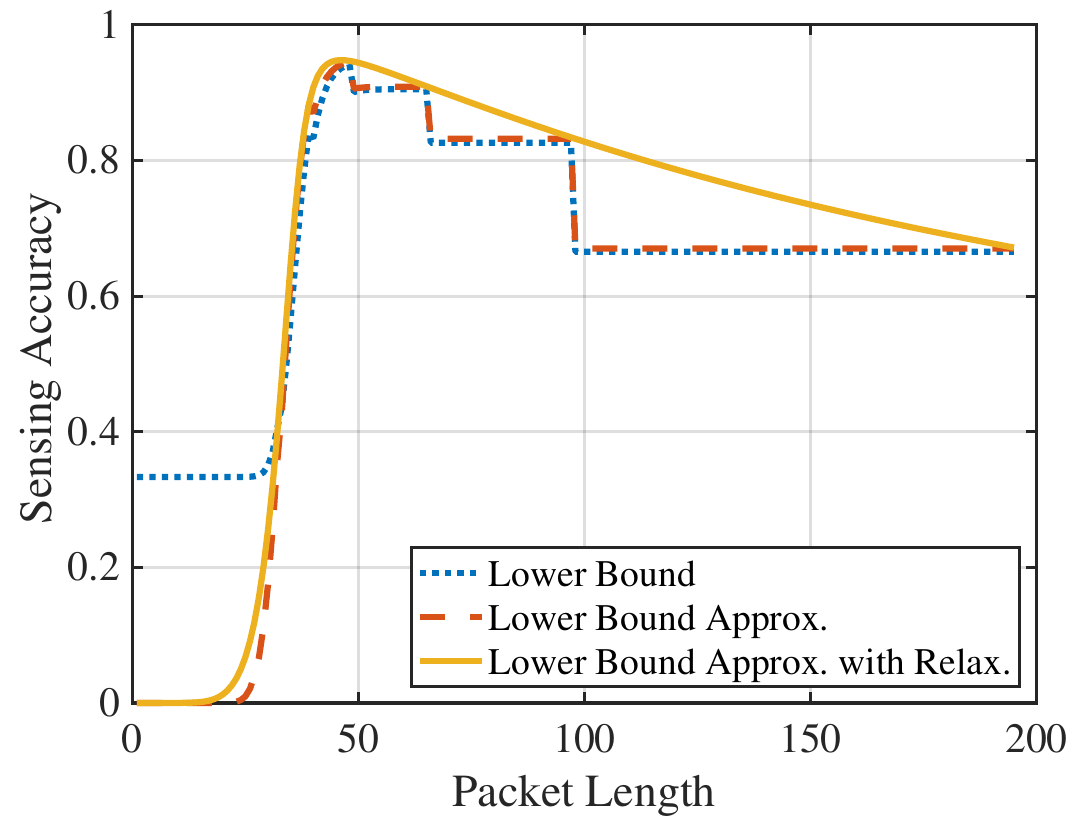}
\label{fig:NumRes_MV_PL_ACC(5dBL3)}
}
\subfigure[3 classes with $\gamma=10\text{dB}$ ]{
\includegraphics[width=0.45\columnwidth]{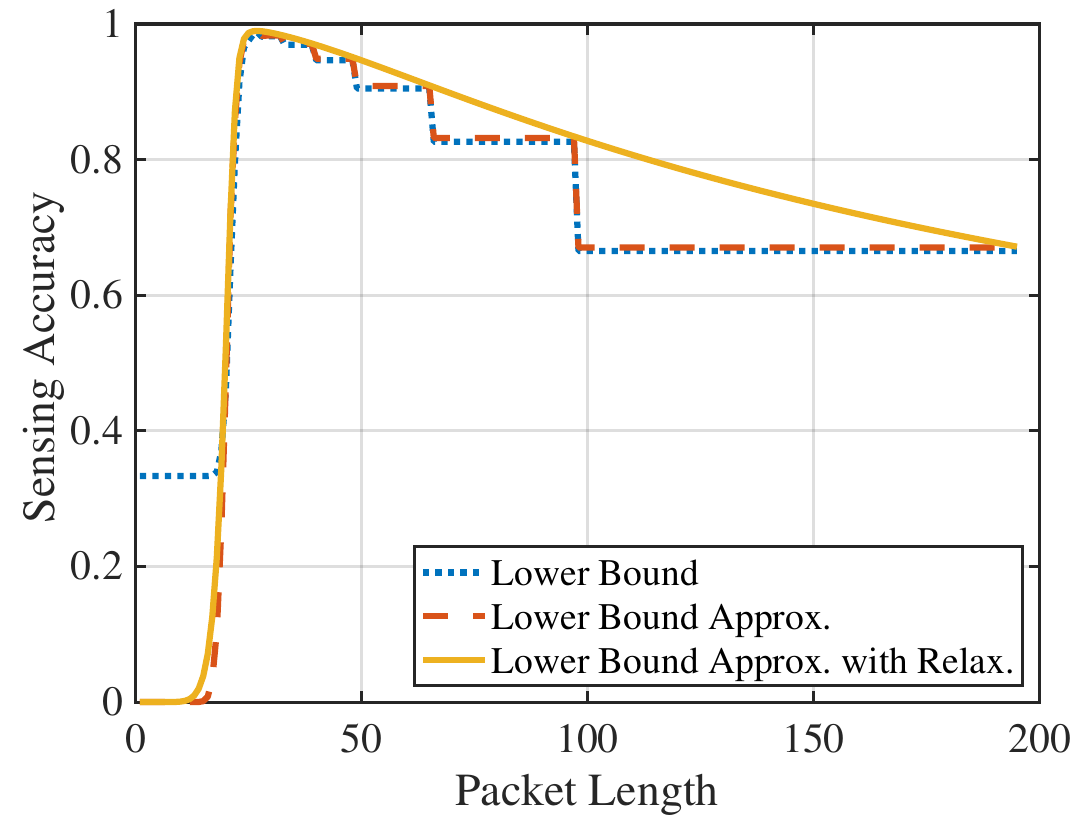}
\label{fig:NumRes_MV_PL_ACC(10dBL3)}}
\caption{Comparison among lower bound of sensing accuracy $\hat{E}_{\sf{mv}}$ in \eqref{sdp+ms}, its approximation in \eqref{MV_UP} and  the approximation with relaxation  in \eqref{E2E_ACC_MV_approx} in case of multi-view sensing.
The settings are $T=1$ms, $\Delta_S=0.02$ms, $B_W=0.2$MHz, $\xi_a=0.95$, $G=1,N=10,Q_B=8,\eta =1.7$.}
\label{Fig:UP_SDP}\vspace{-3mm}
\end{figure}

By relaxing the constraint of $[\cdot]_{(1/L)^+}$ in the lower bound on classification accuracy, we 
obtain the following insightful approximation on \eqref{sdp+ms} that follows directly from the application of first-order Taylor expansion around the mean of the random variable $|\hat{\mathcal{K}}|$:
\begin{equation}
    \mathbb{E}_{|\hat{\mathcal{K}|}}\left[ Q\left(-G\sqrt{|\hat{\mathcal{K}}|}\right)  \right] \approx 
     Q\left(-G\sqrt{\mathbb{E}[|\hat{\mathcal{K}}|]}\right).
\end{equation}
Using this, the lower bound of sensing accuracy for multi-view sensing in \eqref{sdp+ms}, i.e., $\hat{E}_{\sf{mv}}$, can be approximated as
 \begin{equation}
    \label{MV_UP}
    \begin{split}
\hat{E}_{\sf{mv}} \approx   (L-1)Q(-G\sqrt{K_{\sf{mv}}\rho})-(L-2),
    \end{split}
\end{equation}
where $G=\frac{1}{2}\sqrt{g_{\min}}$ and $K_{\sf{mv}}\rho $ is the expected number of successfully transmitted views given a total of $K$ transmissions and successful transmission probability $\rho$.
Fig. \ref{Fig:UP_SDP} compares $\hat{E}_{\sf{mv}}$ in \eqref{sdp+ms}  to its approximation  in \eqref{MV_UP}.
For the binary classifier (see Fig. \ref{fig:NumRes_MV_PL_ACC(5dBL2)} and \ref{fig:NumRes_MV_PL_ACC(10dBL2)}), it can be seen that the approximation error is generally very small. The small error stems from the higher-order moments of $|\hat{\mathcal{K}}|$ and has negligible effects on the packet length optimization.
For the 3-class classifier (see Fig. \ref{fig:NumRes_MV_PL_ACC(5dBL3)} and \ref{fig:NumRes_MV_PL_ACC(10dBL3)}), the approximation error is also small, except for small packet lengths where the relaxation of the constraint of $[\cdot]_{(1/L)^+}$ causes a relatively large error.

\begin{Remark}[Reliability-View Tradeoff]
    The approximation of the lower bound on sensing accuracy in \eqref{MV_UP} demonstrates the tradeoff between communication reliability and the number of transmitted views.
    Specifically, $Q(-G\sqrt{x})$ is a monotonically increasing function of $x$. % Thus, the tradeoff can be revealed by expected successfully transmitted views.
    As shown in Fig. \ref{Fig:UP_SDP}, regardless of the steps resulting from the integer relaxation of $K_{\sf{mv}}$, an increasing packet length causes two opposite effects on sensing accuracy: by lengthening the packets, the total number of views decreases (i.e., reducing the classification accuracy due to fewer sensor views) while the successful decoding probability increases (i.e., improved communication reliability). 
This necessitates the optimization of the packet length, as in the case of multi-snapshot sensing. 
\end{Remark}

\section{Packet-Length Optimization for Ultra-LoLa Inference}
\label{Sec:PL-optimization}
In this section, the packet length that controls the reliability-view tradeoff is optimized to maximize the E2E sensing accuracy by deriving accurate, tractable surrogates of the sensing accuracies.
We treat the two sensing scenarios sequentially.

\subsection{Optimal Packet-Length for Multi-Snapshot Sensing}

We start by optimizing the tradeoff from Proposition \ref{Proposition:E2E_Accuracy_SS} to maximize the lower bound on sensing accuracy of multi-snapshot case. 
First, by expanding $K_{\sf{ms}}$ with the expression in \eqref{eq:num_observations}, the lower bounded sensing accuracy is a function of the packet length $D$, denoted as $\Psi_{\sf{ms}}(D)$,  given as
    \begin{equation}
    \small
\label{eq:E2E_accuracy}
    \begin{split}
   \Psi_{\sf{ms}}(D)    & =\\&(L-1)\rho\left[ \frac{1}{L}- Q\left(G\sqrt{\left\lfloor\frac{T}{\Delta_S}-\frac{D}{B_W \Delta_S}\right\rfloor}\right)  \right]_{(1/L)^+}+\frac{1}{L},
    \end{split}
    \end{equation}
    where $G=\frac{1}{2}\sqrt{g_{\min}}$.
The resulting optimization problem for the $L$-class classifier is
\begin{equation}
\label{prob:E2E_acc_single-sensor}
    \begin{split}
        \max_{D\in \mathbb{Z}_+} & \quad \Psi_{\sf{ms}}(D),
            %\text{s.t.}\quad &  D\in \mathbb{Z}_+, %\{1,2,\dots,D_{\max}\}
    \end{split}
\end{equation}
where $\mathbb{Z}_+$ denotes the set of positive integers. In order to solve the problem, we first approximate the Q-function as \cite{RN329}
\begin{equation}
\label{Approx_Q}
    Q(x) \approx 1-\sigma(\eta x),
\end{equation}
where $\sigma(x)=\frac{1}{1+e^{-x}}$ is the sigmoid function and $\eta>0$ is a fitting coefficient. Throughout the paper, we will use $\eta=1.7$, which generally results in a good approximation for a wide range of $x$~\cite{RN329}.
Next, the integer constraint on $K_{\sf{ms}}$ is relaxed by allowing it to take non-integral values, and finally we note that the function $[\cdot]_{{1/L}^+}$ in \eqref{eq:E2E_accuracy} with the approximation in \eqref{Approx_Q} can be replaced by an upper bound on $D$ given by
\begin{equation}
\label{D_max_upper}
D_{\max} =\left\lfloor \left( T-\frac{\Delta_S(\ln{(L-1)})^2}{G^2\eta^2}\right)B_W
\right\rfloor,
\end{equation}
i.e., we constrain $D\in \{1,2,\ldots,D_{\max}\}$.

Using the relaxations above, we obtain the surrogate  objective of problem \eqref{prob:E2E_acc_single-sensor}, denoted as  $ \tilde{\Psi}_{\sf{ms}}(D)$, given as
\begin{equation}
\label{sdp_discrete_rex}
\tilde{\Psi}_{\sf{ms}}(D)=\xi_a (L-1) \tilde{\nu}_{\sf{ms}}(D)+\frac{1}{L},
\end{equation}
where $ \tilde{\nu}_{\sf{ms}}(D)$ is a positive function of $D\in \{1,\ldots,D_{\max}\}$, given as
\begin{equation}
  \label{nu_ms}
  \tilde{\nu}_{\sf{ms}}(D)=\left( \sigma\left(\psi_I(D)\right)-\frac{L-1}{L}  \right)   \sigma(\psi_T(D)).
\end{equation}
The functions $\psi_I(D)$ and $\psi_T(D)$ are defined as
\begin{equation}
  \label{inner functions}
  \begin{split}
    \psi_I(D) &=  G\eta\sqrt{\frac{T}{\Delta_S} -\frac{D}{\Delta_S B_W} },\\
    \psi_T(D) & = \ln(2)\eta \sqrt{\frac{D}{V}}\left( \log_2 (1+\gamma)-\frac{Q_B N}{D} \right).
  \end{split}
\end{equation}

\begin{figure}[t!]
\centering
\subfigure[Low SNR with $\gamma=0\text{dB}$]{
\includegraphics[width=0.45\columnwidth]{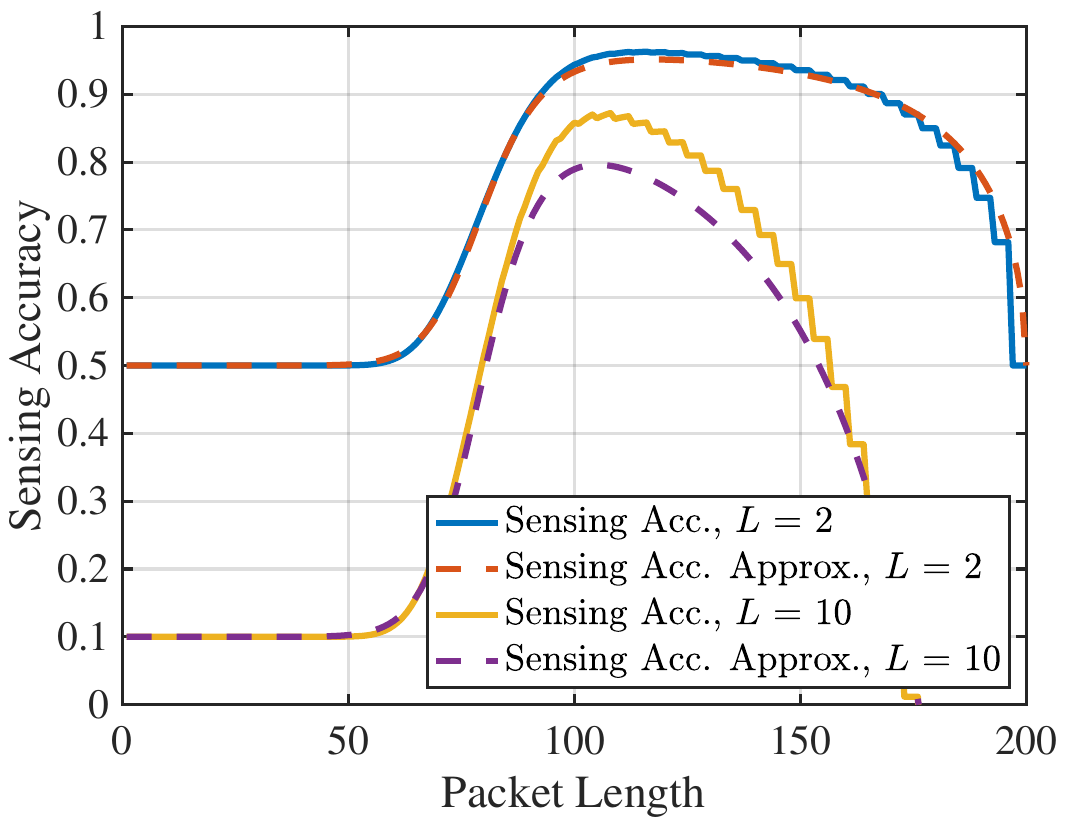}\label{fig:E2E_Acc_Approx_Single_lowsnr}}
\subfigure[High SNR with $\gamma=5\text{dB}$]{
\includegraphics[width=0.45\columnwidth]{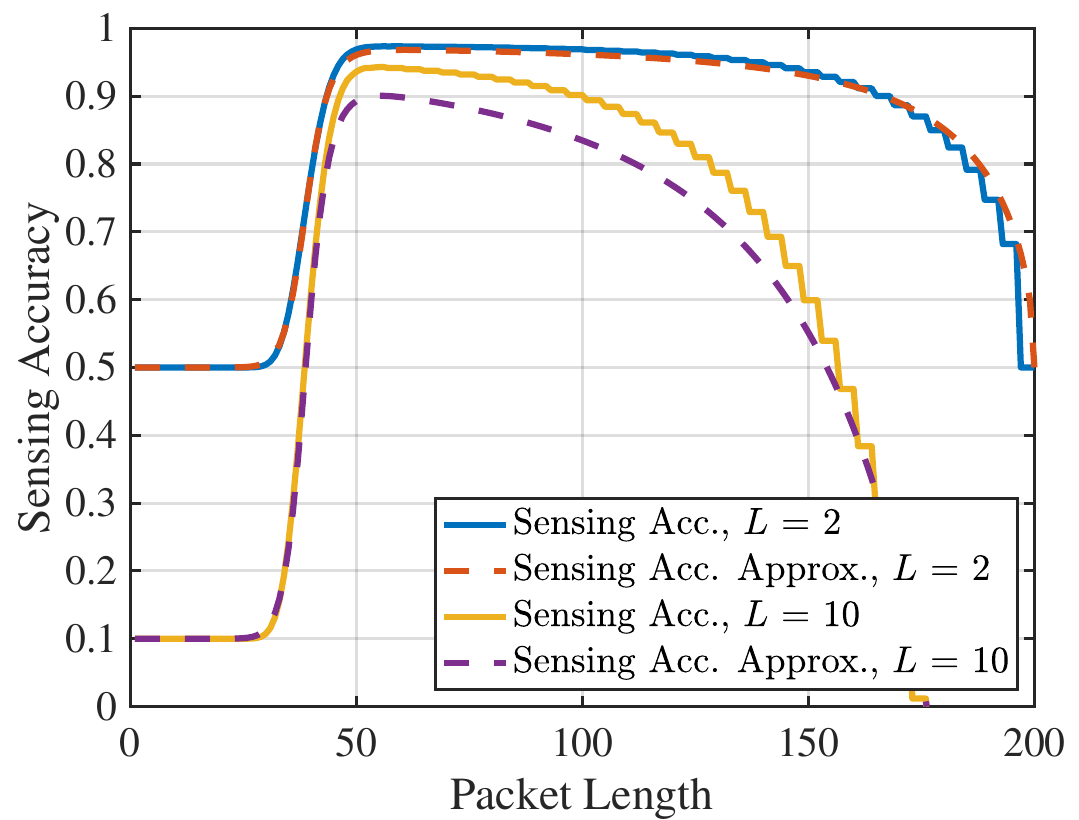}\label{fig:E2E_Acc_Approx_Single_highsnr}}
\caption{Comparison between  \eqref{eq:E2E_accuracy} and its approximation  \eqref{sdp_discrete_rex} in case of multi-snapshot sensing. 
The settings are identical to that of multi-view case in Fig. \ref{Fig:UP_SDP}  except $G=0.5$.
}
\label{fig:E2E_Acc_Approx_Single}
\vspace{-3mm}
\end{figure}

As shown in Fig. \ref{fig:E2E_Acc_Approx_Single}, the proposed function $\tilde{\Psi}_{\sf{ms}}(D)$ provides a good approximation of the lower bound of sensing accuracy $\Psi_{\sf{ms}}(D)$. 
Building on this approximation, we seek the value of $D$ that maximizes $\tilde{\Psi}_{\sf{ms}}(D)$. Since $\tilde{\Psi}_{\sf{ms}}(D)$ is a monotonically increasing function of $\tilde{\nu}_{\sf{ms}}(D)$, the optimal value of $D$ can be found by solving the simplified problem:
\begin{equation}
\label{prob:max_E2E_acc_MC}
    \begin{split}
        \max_D & \quad \tilde{\nu}_{\sf{ms}}(D)\\
            \text{s.t.}\quad &  D\in \{1,2,\dots,D_{\max}\}.
    \end{split}
\end{equation}
$\tilde{\nu}_{\sf{ms}}(D)$ in problem \eqref{prob:max_E2E_acc_MC} is found to be a log-concave function of $D$ with a unique maximum elaborated in Proposition \ref{Optimal_BL_Single_sensor}. 

% if and only if
% \begin{equation}
% \label{condtion:exist_optimalBL_singlesensor}
%    f_{\bf{ms}}(1)\cdot f_{\bf{ms}}(D_{\max} )<0,
% \end{equation}
% where $f_{\bf{ms}}(D)$ is a function of $D$ determining the monotone of the objective function, given as
% \begin{equation}
%        f_{\bf{ms}}(D)=\frac{\psi'_I(D)}{(1+e^{\psi_I(D)})(1/L-e^{-\psi_I(D)})}+\frac{\psi'_T(D)}{1+e^{\psi_T(D)}}.
% \end{equation}
% with the  first derivatives:   $\psi'_I(D)=-\frac{G\eta  }{2\sqrt{\Delta_SB_W(B_WT-D)}}<0, \psi'_T(D)=\frac{ \ln{(2)}\eta}{2\sqrt{DV}}\left(\log_2(1+\gamma)+\frac{Q_BN}{D}  \right)>0$.

\begin{Proposition}[Optimal Packet Length for Multi-Snapshot Sensing]
\label{Optimal_BL_Single_sensor}
Let
\begin{equation}
\begin{split}
\small
       f_{\bf{ms}}(D)=&C_{\bf{ms}}\sqrt{\frac{B_WT}{D}-1}\left( \log_2(1+\gamma)+\frac{Q_BN}{D}\right)\\ 
       &- \frac{1+e^{\psi_T(D)}}{(1+e^{\psi_I(D)})(1/L-e^{-\psi_I(D)})},
       \end{split}
\end{equation}
where $C_{\bf{ms}}=\frac{\ln(2)}{G}\sqrt{\frac{\Delta_S B_W}{V}}$ is a positive constant and $\psi_T(D),\psi_I(D)$ are defined in \eqref{inner functions}.
The optimal packet length that solves problem \eqref{prob:max_E2E_acc_MC} is then
\begin{equation}
   D^*= 
   \left \lfloor \tilde{D}^* \right\rceil_{ \tilde{\nu}_{\sf{ms}}(\cdot)},
\end{equation}
where $\lfloor x \rceil_{\tilde{\nu}_{\sf{ms}}(\cdot)}$ is equal to $\lfloor x \rfloor$ if $ \tilde{\nu}_{\sf{ms}}(\lfloor x \rfloor)\geq \tilde{\nu}_{\sf{ms}}(\lceil x \rceil) $, and is otherwise equal to $\lceil x \rceil$, and $\tilde{D}^*$ is given as
\begin{equation}
 \tilde{D}^*  = \left\{D|  f_{\bf{ms}}(D)=0,D\in [1,D_{\max}]  \right\}
\end{equation}
if $f_{\bf{ms}}(1)\cdot f_{\bf{ms}}(D_{\max} )<0$ holds, otherwise $\tilde{D}^*  = \argmax_{D\in \{1,D_{\max} \}} \tilde{\nu}_{\sf{ms}}(D)$.

% \begin{equation}
%    \tilde{D}^*= \begin{cases}
%    \left\{D|  f_{\bf{ms}}(D)=0,D\in [1,D_{\max}]  \right\},& \text{if } f_{\bf{ms}}(1)\cdot f_{\bf{ms}}(D_{\max} )<0,\\
%    \argmax_{D\in \{1,D_{\max} \}} \tilde{\nu}_{\sf{ms}}(D) , &\text{otherwise}.
%     \end{cases}
% \end{equation}
\end{Proposition}
\noindent The proof is given in Appendix \ref{Proof_Existence_Optimal_BL}.

% \begin{Remark}[Optimal packet Length of Binary Classifier]
    
% \end{Remark}

\begin{figure}[t!]
\centering
\subfigure[Optimal packet length vs. SNR ]{
\includegraphics[width=0.45\columnwidth]{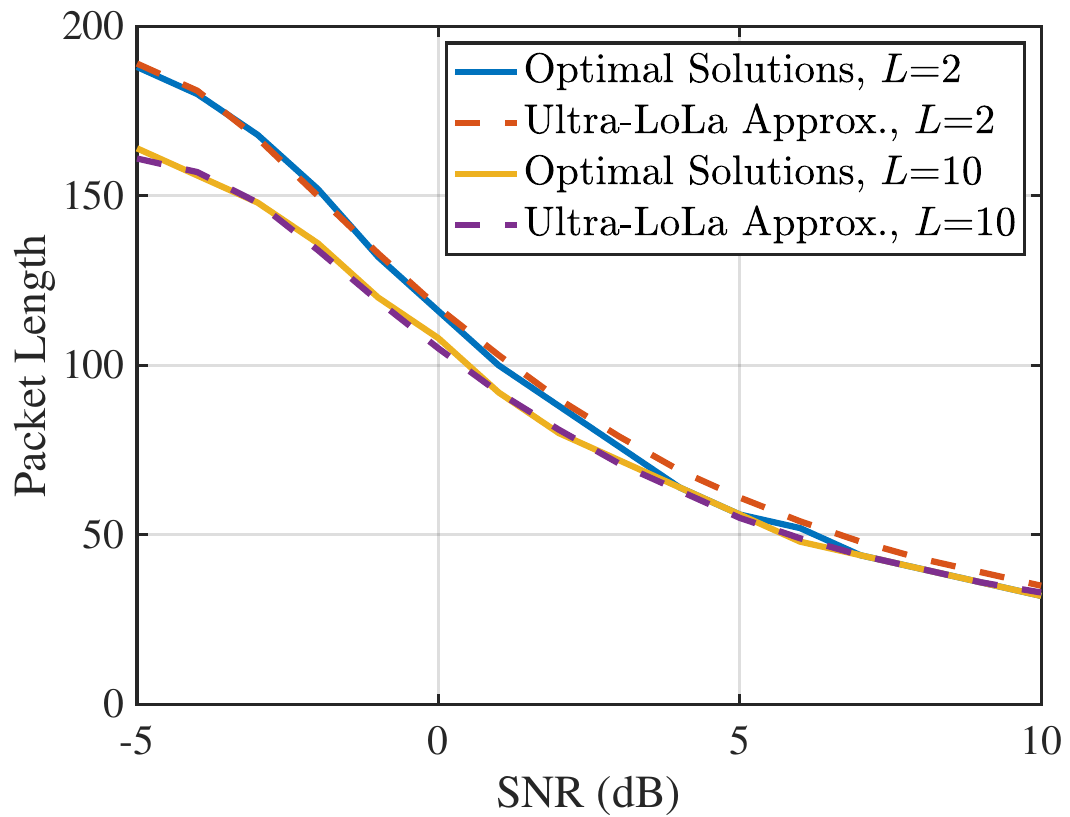}\label{Fig: SNR_OpBL_exp}}
\subfigure[Sensing accuracy vs. SNR]{
\includegraphics[width=0.45\columnwidth]{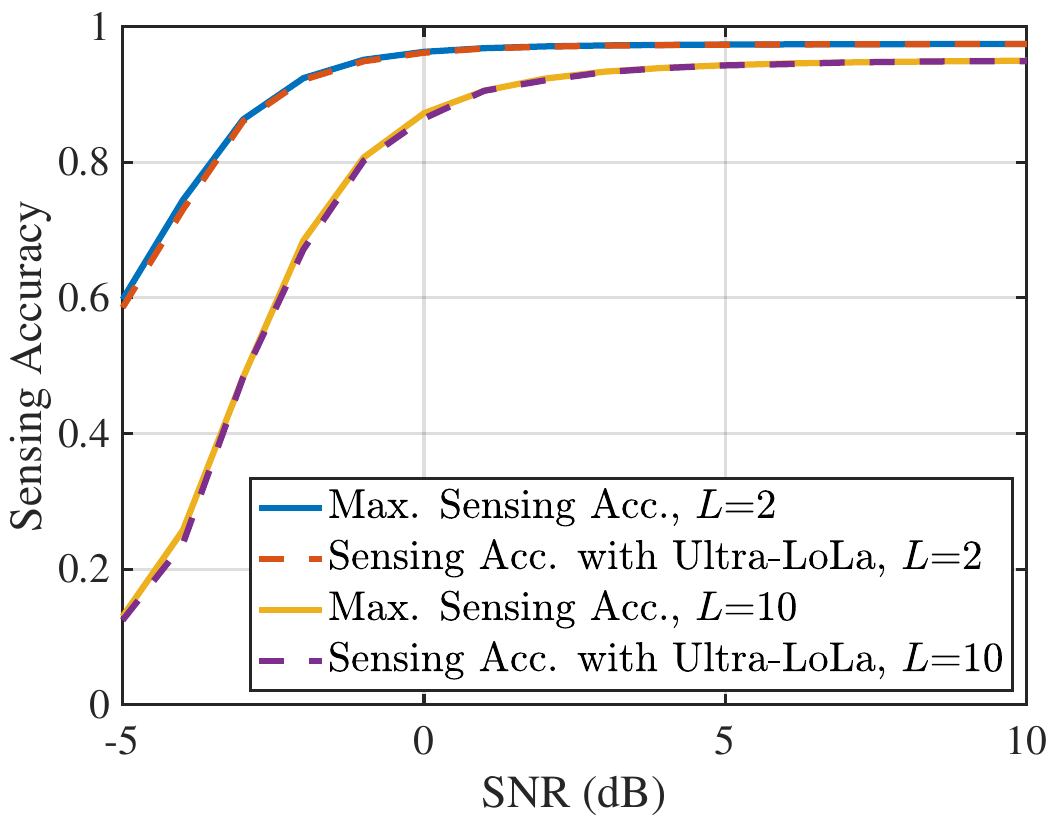}\label{Fig: SNR_SDP_exp}}
\caption{Comparison between optimal solutions and approximated solutions using ultra-LoLa scheme 
in the case of multi-snapshot sensing
with the same setting as Fig. \ref{fig:E2E_Acc_Approx_Single}. 
}
\label{Fig:SNR effects}\vspace{-4mm}
\end{figure}

Fig. \ref{Fig:SNR effects} compares the ultra-LoLa packet length obtained using Proposition \ref{Optimal_BL_Single_sensor} to the optimal solution of $\Psi_{\sf{ms}}(D)$.
The approximations are close to the optimal ones, which validates the surrogate objective function as an accurate approximation.
The figure also shows that the optimal balance between communication reliability and classification accuracy is controlled by the SNR and the number of classes $L$.
On the one hand, the optimal packet length decreases as the SNR increases, allowing for more snapshots and a higher  classification accuracy. On the other hand, increasing the number of classes decreases the optimal packet length as more snapshots are required for a classification task with more classes.

\subsection{Optimal Packet-Length for Multi-View Sensing}

We now turn to the problem of choosing the optimal packet length of multi-view sensing. 
First, we expand $K_{\sf{mv}}$ in \eqref{MV_UP} with expression in \eqref{MS_K} and formulate the optimization problem of interest, given as
\begin{equation}
\label{PF_MS_SDP_UP}
    \begin{split}
        \max_{D\in\mathbb{Z}_+} &  \quad \Psi_{\sf{mv}}(D),%\\
            %\text{s.t.}\quad &  D\in \mathbb{Z}_{+}.%\{1,\dots,D_{\max}\},
    \end{split}
\end{equation}
where $ \Psi_{\sf{mv}}(D)$ is the approximation of the lower bounded sensing accuracy:
\begin{equation}
   \Psi_{\sf{mv}}(D) =   (L-1)Q(-G\sqrt{\nu_{\sf{mv}}(D)})-(L-2),
\end{equation}
with $\nu_{\sf{mv}}(D)$ being the expected number of successfully transmitted views, given as
\begin{equation}
\label{Exp_SenNum}
     \nu_{\sf{mv}}(D)=\left\lfloor \frac{(T-\Delta_S)B_W}{D} \right\rfloor\rho.
\end{equation}
Next, we apply again the approximation of Gaussian Q-function in \eqref{Approx_Q} and a relaxation of an integral number of sensors. 
With the relaxation and approximation above, we obtain the surrogate objective function of problem \eqref{PF_MS_SDP_UP}, denoted as $\Tilde{\Psi}_{\sf{mv}} (D) $, given as
\begin{equation}
\label{E2E_ACC_MV_approx}
\Tilde{\Psi}_{\sf{mv}} (D) = (L-1)Q(-G\sqrt{\Tilde{\nu}_{\sf{mv}}(D)})-(L-2),
\end{equation}
where $ \Tilde{\nu}_{\sf{mv}}(D)$ is the approximation of the expected number of received views, given as
\begin{equation}
    \Tilde{\nu}_{\sf{mv}}(D)=\frac{(T-\Delta_S)B_W\xi_a}{D\left(1+\exp(-\psi_T(D))\right)},
\end{equation}
with $\psi_T(D)$ defined in \eqref{inner functions}.

The approximation error of the approximate sensing accuracy is illustrated in Fig. \ref{Fig:UP_SDP}.
It is seen that the continuous relaxation in \eqref{E2E_ACC_MV_approx} captures the monotonicity and approximates the optimum of the exact expression.
Using this approximation, we optimize the packet length $D$ to maximize  $\Tilde{\Psi}_{\sf{mv}}(D)$. Since $Q(-G\sqrt{x})$ is a monotonically increasing function of $x$, the optimization problem on $\Tilde{\Psi}_{\sf{mv}}(D)$ reduces to 
\begin{equation}
\label{prob:max_E2E_acc_MV}
    \begin{split}
        \max_D & \quad \tilde{\nu}_{\sf{mv}}(D)\\
            \text{s.t.}\quad &  D\in \{1,2,\dots,D_{\max}\},
    \end{split}
\end{equation}
where $D_{\max}=\left\lfloor (T-\Delta_S)B_W  \right\rfloor$ is the maximum packet length ensuring $K\geq 1$ (clearly, $K=0$ can never be better than $K\ge 1$).

% We arrive to the following optimization problem:
% \begin{equation}
% \label{PF_MS_Exp_SenNum}
%     \begin{split}
%         \max_D &  \quad \Tilde{\Psi}_{\sf{mv}} (D)\\
%             \text{s.t.}\quad &  D\in \{1,\dots,D_{\max}\},
%     \end{split}
% \end{equation}
% where
% \begin{equation}
% \label{E2E_ACC_MV_approx}
% \Tilde{\Psi}_{\sf{mv}} (D) = (L-1)Q(-G\sqrt{\Tilde{\nu}_{\sf{mv}}(D)})-(L-2)
% \end{equation}
% with $\Tilde{\nu}(D)$ expressed as
% \begin{equation}
%     \Tilde{\nu}_{\sf{mv}}(D)=\frac{(T-\Delta_S)B_W}{D\left(1+\exp(-\psi_T(D))\right)},
% \end{equation}
% where $\psi_T(D)$ is defined in \eqref{inner functions}.

Finally, $\tilde{\nu}_{\sf{mv}} (D)$ is found to be an unimodal function of $D\in [1,D_{\max}]$ with a unique maximum.
Considering the feasible packet length $ D \in \{1,\ldots, D_{\max} \}$, the optimal packet length of problem \eqref{prob:max_E2E_acc_MV}, denoted as $D^*$, is in the interior of the interval if and only if the derivative of $\tilde{\nu}_{\sf{mv}} (D)$ at the two boundary points have different signs.
This is satisfied when
\begin{equation}
  \label{condtion_optimalBL_MS}
  f_{\bf{mv}}(1)\cdot f_{\bf{mv}}(D_{\max} )<0,
\end{equation}
where $ f_{\bf{mv}}(D)$ is a function of $D$, given as
\begin{equation}
  f_{\bf{mv}}(D)=\frac{\eta \ln(2) \left( \log_2(1+\gamma)+NQ_B/D \right)}{2\sqrt{DV}(\exp(\psi_T(D))+1)}-\frac{1}{D}.
\end{equation}
Otherwise, the optimal packet length is  obtained at the endpoints, i.e., 
\begin{equation}
  D^*={\argmax}_{D\in \{1,D_{\max} \}  } \tilde{\nu}_{\sf{mv}}(D)  .
\end{equation}
The optimal packet length given that condition \eqref{condtion_optimalBL_MS} is satisfied is presented in Proposition \ref{Unique_Solution_MS}.
  
\begin{Proposition}[Optimal Packet Length for Multi-View Sensing]
\label{Unique_Solution_MS}
If the condition in \eqref{condtion_optimalBL_MS} is satisfied, the optimal packet length $D^*$ can be obtained as
\begin{equation}
   D^*= 
   \left \lfloor\frac{ V \zeta^2 }{\ln^2(2)\eta^2\log^2_2(1+\gamma)}\right\rceil_{ \Tilde{\nu}_{\sf{mv}}(\cdot)},
\end{equation}
where $\lfloor x \rceil_{\tilde{\nu}_{\sf{mv}}(\cdot)}$ is equal to $\lfloor x \rfloor$ if $ \tilde{\nu}_{\sf{mv}}(\lfloor x \rfloor)\geq \tilde{\nu}_{\sf{mv}}(\lceil x \rceil) $, and is otherwise equal to $\lceil x \rceil$, and $\zeta$ is the solution to the transcendental equation
\begin{equation}
\label{transcendental}
   \zeta=\left\{ x|x+\frac{\omega}{x}-2\exp{\left( x- \frac{\omega}{x}\right)}=2 ,x>0\right \},
\end{equation}
with $\omega=\frac{\ln^2(2)\eta^2}{V}\log^2_2(1+\gamma) N  Q_B$ being a constant representing the effects of channel reliability and sensing quality.
    
\end{Proposition}
\noindent The proof is given in Appendix \ref{Proof_Unique_Solution_Prob}.
\begin{figure}[t!]
\centering
\subfigure[ Optimal packet length vs. SNR ]{
\includegraphics[width=0.45\columnwidth]{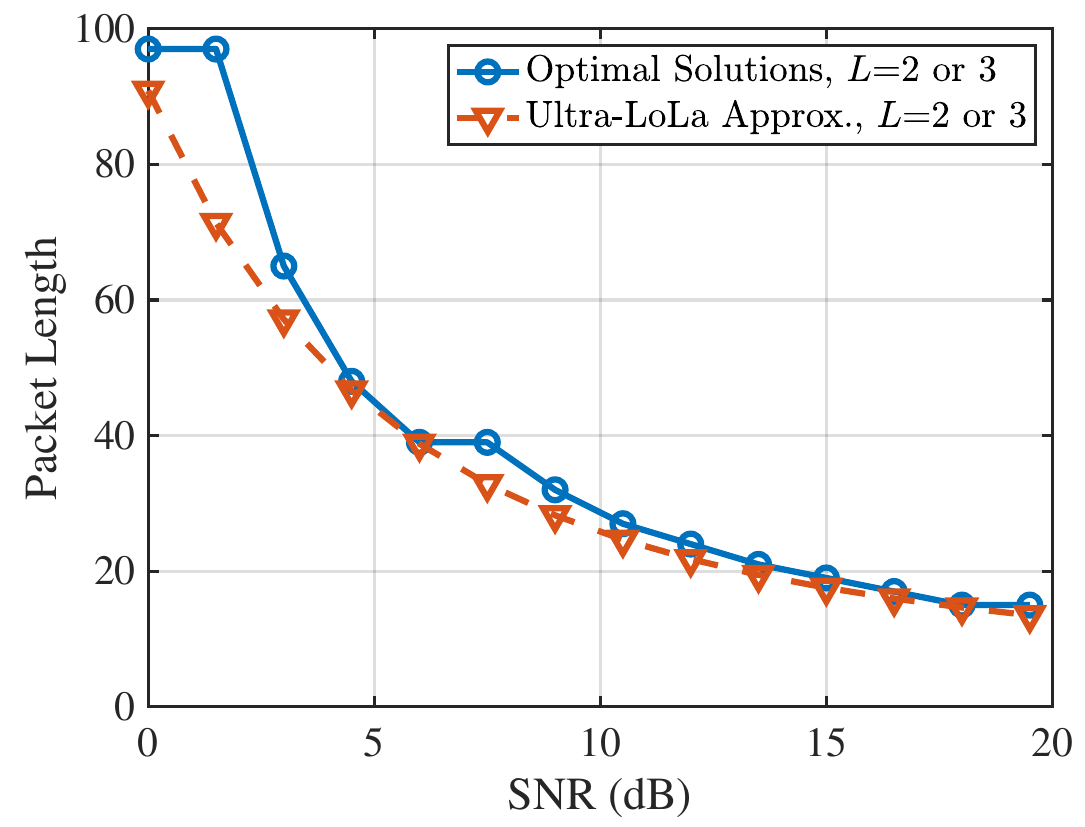}\label{Fig: SNR_OpBL_MS}}
\subfigure[ Sensing accuracy vs. SNR]{
\includegraphics[width=0.45\columnwidth]{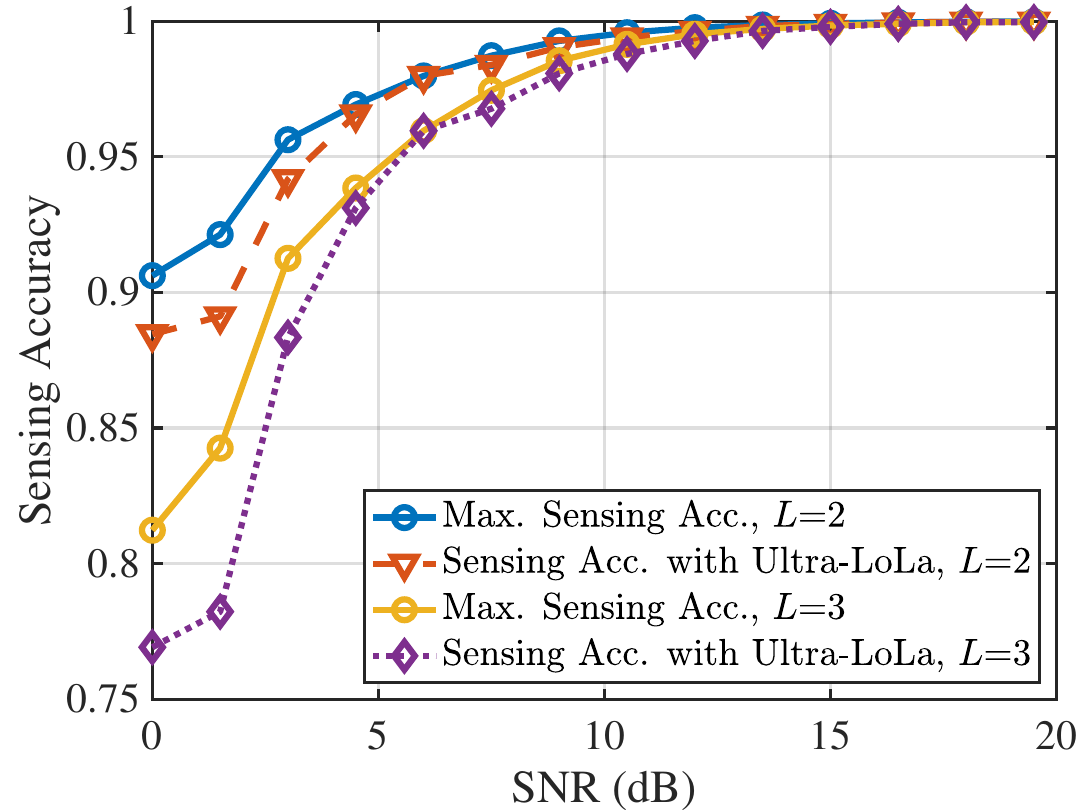}\label{Fig: SNR_SDP_MS}}
\caption{Comparison between optimal solutions and approximate solutions using ultra-LoLa scheme in case of multi-view sensing and  the same settings as Fig. \ref{Fig:UP_SDP}.}
\label{Fig:SNR effects_MS}
\vspace{-4mm}
\end{figure}

Fig. \ref{Fig:SNR effects_MS} compares the ultra-LoLa packet length obtained by solving the equation in Proposition \ref{Unique_Solution_MS} to the optimal solution in the settings of SNR and number of classes.
Specifically, Fig. \ref{Fig: SNR_OpBL_MS} demonstrates the trend of decreasing optimal packet length with the growth of SNR. 
This is because that the improvement in packet decoding from increasing SNR can allow shorter packets  such that more sensors can be scheduled before the end of the task.
Otherwise, the approximate solution is seen to be smaller than the optimal counterpart. The reason is that the optimal packet length guaranteeing an integer total number of sensors is larger than its relaxation. 
% By increasing SNR, the gap between optimal and approximated solution is reduced. 
Additionally, both 3-class and binary classification tasks  have the same solutions since Problem \eqref{prob:max_E2E_acc_MV} targets maximizing the expected number of active sensors $\Tilde{\nu}_{\sf{mv}}(D)$, which is independent of the number of classes.
Correspondingly, in Fig. \ref{Fig: SNR_SDP_MS}, the sensing accuracy with the approximate solution is observed to approach its optimal counterpart with a shrinking gap as the SNR increases.
Overall, both an increase in SNR and a decrease in the number of classes enhance the E2E performance due to higher successful decoding probability and an easier classification task, respectively.

Determining the optimal packet length for multi-snapshot and multi-view sensing can be approximated by maximizing the associated surrogate functions via bisection search over the feasible set $D\in \{1,2,\dots, D_{\max}\}$. The resulting computational complexity is 
$\mathcal{O}(\log D_{\max})$, which is substantially lower than the 
$\mathcal{O}\{D_{\max}K_{\max}\}$ complexity of brute-force search, where $K_{\max}$ denotes the maximum number of views.
In the context of practical realization in 5G systems, the optimized packet length computed by server can be assigned to sensors using the \emph{physical downlink control channel} (PDCCH) signaling \cite{3GPP_38_212}.

\ifthenelse{\boolean{showaircomp}}{
\section{Ultra-LoLa Inference with Over-the-Air Computing}
\label{AirComp}

Different from the orthogonal access in Section \ref{Sec:Multiview}, this section considers simultaneous feature transmission over distributed sensors using AirComp. 
To realize the fair comparison between orthogonal access and AirComp in the E2E latency, we consider a simple method of repetition coding to improve the reliability of AirComp via repeatedly transmitting feature vector \cite{liu2024digital}.
In particular, sensor $k$ transmits $\mathbf{x}_k$ with repetition of $R$ times using linear modulation and full bandwidth $B_W$.
Given task completion deadline $T$, the repetition time $R$ is maximized as  $R=\lfloor \frac{(T-\Delta_S)B_W}{N} \rfloor$.
Building on the setting above, the analysis results are provided as follows.

Considering the received feature vector in \eqref{AirComp features} and repetition coding, the covariance matrix of $\overline{\mathbf{x}}$, denoted as $\widetilde{\mathbf{C}}$, is given as 
\begin{equation}
\widetilde{\mathbf{C}}=\frac{\mathbf{C}}{|{\mathcal{K}}|}+ \frac{1}{2 |{\mathcal{K}}|^2R\gamma'} \mathbf{I}. 
\end{equation}
Feeding the noisy feature vector $\overline{\mathbf{x}}$ into the classifier, the resulting minimum discrimination gain for this case, denoted as $\Tilde{g}_{\min}$,  can be expressed as
\begin{equation}
\label{eq:scaled_DG_AirComp} \Tilde{g}_{\min}=\min_{\ell,\ell'}\{ [\boldsymbol{\mu}_{\ell}-\boldsymbol{\mu}_{\ell'}]^{\sf{T}}\widetilde{\mathbf{C}}^{-1}[\boldsymbol{\mu}_{\ell}-\boldsymbol{\mu}_{\ell'}]\},
\end{equation}
Substituting  \eqref{eq:scaled_DG_AirComp} into \eqref{sdp+ms}, the approximation of expected accuracy for AirComp case, denoted as $\overline{E}_{\sf{air}}$, is computed as
\begin{equation}
\label{Aircomp_E2E}
    \begin{split}
  \overline{E}_{\sf{air}}&= (L-1) 
          \sum_{|{\mathcal{K}}|=0}^K Q\left( -\frac{\sqrt{ \Tilde{g}_{\min}} }{2} \right) {P}_{|{\mathcal{K}}|} -(L-2),
    \end{split}
\end{equation}
where ${P}_{|{\mathcal{K}}|}=\binom{K}{|{\mathcal{K}}|} \xi_a ^{|{\mathcal{K}}|}(1-\xi_a)^{K-|{\mathcal{K}}|}$  is the PMF of the binomial distribution $|{\mathcal{K}}|\sim  \mathcal{B}(K,\xi_a)$.

To visualize the effect of sensor number on  inference accuracy, we consider a special case of classifier and provide the sensing accuracy of AirComp case in Lemma \ref{Aircomp_E2E_SP}.

\begin{Lemma}[Sensing Accuracy of Multi-view with AirComp]
\label{Aircomp_E2E_SP}
Considering a classifier with covariance matrix $\mathbf{C}=\sigma^2_c\mathbf{I}$, the sensing accuracy in \eqref{Aircomp_E2E} can be approximated as
\begin{equation}
\label{eq:analog_E2E_acc} \overline{E}_{\sf{air}}= (L-1) 
           Q\left( -\frac{1}{2}\sqrt{\frac{| \mathcal{K}|^2\boldsymbol{\mu}_{\min}}{|\mathcal{K}|\sigma_c^2+\frac{1}{2 R\gamma'}}}  \right)  -(L-2),
\end{equation}
where $\boldsymbol{\mu}_{\min}=\min_{\ell,\ell'}\{\Vert\boldsymbol{\mu}_{\ell}-\boldsymbol{\mu}_{\ell'}\Vert^2\}$ is the minimum distance between ground-truths of arbitrary two classes;$|\mathcal{K}|$ is the number of activation sensors.
\end{Lemma}

With the same setting as Lemma \ref{Aircomp_E2E_SP}, for the case of $K$ sensors performing short-packet transmission, the  sensing accuracy of multi-view using ultra-LoLa packet length, can be approximated as
\begin{equation}
\label{Digital_IG}
\begin{split} \overline{E}_{\sf{mv}}\approx (L-1)Q\left(-\frac{1}{2} \sqrt{ \frac{|\mathcal{K}|\sigma(\psi_T(D^*))\boldsymbol{\mu}_{\min}}{\sigma_c^2}  } \right)-(L-2),
\end{split} 
\end{equation}
where $D^*$ is the packet length obtained in Proposition \ref{Unique_Solution_MS}.

In \eqref{eq:analog_E2E_acc} and \eqref{Digital_IG}, increasing number of devices can both be seen a rising trend of sensing accuracy, while channel noise and decoding error probability degrade the performance of analog and digital transmission, respectively. 
On the other,  repetition in AirComp (i.e., longer deadline or higher bandwidth) leads to an accuracy improvement.

Then we compare the sensing accuracy of Ultra-LoLa with its AirComp counterpart under the power budget $P_{\max}$ and E2E latency constraint $T$ when considering AWGN channel (i.e., $\xi_a=1$). 
The condition of preference for AirComp is provided in Proposition \ref{AnalogVsDigital}.
\begin{Theorem}[When AirComp is Better?] 
    \label{AnalogVsDigital}
 Given the AWGN channel with  power budget $P_{\max}$ and the task completion deadline $T$,  AirComp outperforms  ultra-LoLa short packet transmission if and only if the ultra-LoLa packet length $D^*$ subjects to
    \begin{equation}
    \label{ratio_comparison}
      D^* \leq \frac{1}{4} \left({\beta}+\sqrt{\beta^2+4\lambda}\right)^2,
    \end{equation}
where $\beta=\frac{\sqrt{V}\ln(2\sigma_c^2R\gamma)}{\eta\ln(1+\gamma)}$ and $\lambda=\frac{Q_BN}{\log_2(1+\gamma)}$ with  $\gamma=\frac{P_{\max}}{N_0B_W}$ being the receive SNR of the two cases and $R=\left\lfloor\frac{(T-\Delta_S)B_W}{N}\right\rfloor$ being maximum repetition time ensuring the same E2E latency. 
\end{Theorem}
\noindent The proof is given in Appendix \ref{Proof_AnalogVsDigital}.

Considering the fair constraints in power budget and latency constraint, Proposition \ref{AnalogVsDigital} 
suggests a threshold-based selection scheme between AirComp and ultra-LoLa SPT.
Specifically, when the packet length of ultra-LoLa scheme below the threshold, AirComp is preferred. 
This is because shorter packet results in higher decoding error, thus reducing the sensing accuracy, while AirComp using analog transmission would not be affected.
Among others, if a high quantization resolution of feature vector is required, the AirComp scheme attains better performance due to the analog transmission that does not require the bit quantization.
}

% Given the parameters of orthogonal access, like SNR $\gamma$, packet length $D$ and feature size $N$, if the receive SNR of AirComp case $\gamma'$ is larger than the expression of orthogonal access $\frac{D\exp(\psi_T(D))}{2(T-\Delta_S)B_W\sigma_c^2}$, AirComp is preferred. 
% Otherwise, orthogonal access obtains better performance.
% The reason is that, compared with decoding error in digital mode, the
% effects of the induced channel noise on sensing accuracy is  mitigated by increasing receive SNR of AirComp transceiver.
% Others, the threshold is seen as a monotone increasing function of packet length $D$. This implies that short-packet transmission outperforms AirComp when the packet length is large enough. This comes from the enhancement of successful decoding probability of lengthened packet. 

\section{Experimental Results}
\label{Experimental_res}

In this section, we demonstrate the advantages of the proposed framework through numerical results. We first outline the setup, followed by results for linear classification with the GMM. Finally, we demonstrate how the framework performs in MVCNN-based classification.

\subsection{Experimental Setup}
Unless specified otherwise, the experimental settings are set as follows.
\subsubsection{System and Communication Settings}
We consider a distributed sensing system comprising a single server. For multi-snapshot sensing, we consider a single sensor, while for multi-view sensing, we consider $K$ distributed sensors.
Constrained by limited bandwidth due to the coexistence of other services, the sensing task operates over a bandwidth of $B_W = 100$~kHz for packet transmission, utilizing the access protocols described in Section~\ref{sec:access_protocols}.
For each sensor, we fix the activation probability to $\xi_a=1-10^{-5}$, guaranteeing a low outage probability as required by URLLC.
The packet decoding error probability is assumed to be determined by \eqref{eq:varepsilon}.

\subsubsection{Sensing  and Classification Settings}

We consider a very stringent URLLC task completion deadline of $T = 1$ ms encompassing both sensing and feature transmission phases.
The sensing time consists of camera sensing and on-sensor feature extraction is set as $\Delta_S = 0.1$ ms (representing, for instance, an industrial-level camera \cite{Weinberg:20} and lightweight local computation~\cite{NVIDIA_Jetson_TX2}).

%Sensor scheduling is assumed to be error-free and instantaneous. The relation between the number of images and packet length is described in \eqref{eq:num_observations} and \eqref{MS_K}, for multi-snapshot and multi-view sensing, respectively.
To meet the strict deadline, we assume a quantization resolution of 4 bits per feature entry, which has proven to be sufficient for the considered inference tasks and is also consistent with the literature \cite{sun2020ultra}.
The linear and MVCNN classifiers are specified below.
\begin{itemize}
    \item \textbf{Linear Classification on Synthetic GMM Data:} For the linear GMM classifier, we extract feature vectors according to \eqref{data_dis} with $L=2,5,10$ classes (selected to vary the classification complexity) and the feature dimensionality set to $N=20$. The centroids of the $L$ clusters are specified as follows: for class $\ell$ ($\ell=1,\ldots,L$), the elements from dimension $N(\ell-1)/L$ to dimension $N\ell/L$ are set to $-1$, while the remaining dimensions are set to $+1$. The covariance matrix is defined as $\mathbf{C}=3\mathbf{I}_{N}$.
%The inferred label is obtained using \eqref{Maha_min_classifier}.

\item \textbf{Non-Linear MVCNN-Based Classification on Real-World Data:}
For the MVCNN classifier, we consider the well-known ModelNet dataset \cite{ModelNet-Ref}, which comprises multi-view images of objects (e.g., a person or a plant), and the popular VGG16 model \cite{simonyan2015deep} to implement the MVCNN architecture.
We adapt the VGG16 model for our purposes by separating it into feature extractor and classifier networks, with the classifier running on the server and the feature extractor running on the sensors as in \cite{Zhiyan-AirPooling}. 
The resulting MVCNN architecture is trained for average pooling and targets a subset of ModelNet, comprising $L=20$ popular object classes.
To simulate the sensing phase, each of the $K$ sensors randomly draws an image (view) of the same class (selected uniformly at random) from the dataset.
To reduce the communication overhead, each ModelNet image is resized from $3\times224 \times224$ to  $3\times56 \times56$ (i.e., reducing the image resolution)  and fed into the on-device feature extractor, resulting in a $512\times 1 \times 1$  tensor. This tensor is then further compressed by a fully connected layer to obtain an output feature vector of dimension $N=20$.
\end{itemize}

\subsubsection{Benchmarks}
We consider the following benchmarks:
\begin{itemize}
 \item \textbf{Brute-force Search:} The optimal packet length is obtained by an exhaustive search among all feasible solutions to guarantee the maximum sensing accuracy.
    \item  \textbf{URLLC:} The packet length is determined by computing the shortest packet meeting the decoding error threshold of  $10^{-5}$ commonly required in URLLC \cite{Petar-IEEEProc-2016}.
    \item \textbf{Shannon Rate:} The packet length, denoted as $D_{\sf{sr}}$, is determined based on the Shannon rate at the target SNR, i.e.,
    $D_{\sf{sr}}=\left\lceil \frac{NQ_B}{\log_2(1+\gamma)} \right\rceil$. 
    % Given the transmission duration $\frac{D_{\sf{sr}}}{B_W}$ of each packet, the number of views is maximized to meet the task completion time $T$.
\end{itemize}
Given the task deadline of $T=1$~ms and pre-defined packet lengths in the URLLC and Shannon rate baselines, the number of views for multi-snapshot and multi-view sensing is maximized with respect to \eqref{eq:num_observations} and \eqref{MS_K}, respectively.

\subsection{Ultra-LoLa Inference with Linear Classification}

\begin{figure}[t!]
\centering
\subfigure[Multi-snapshot sensing (SNR = 5dB)]{
\includegraphics[width=0.45\columnwidth]{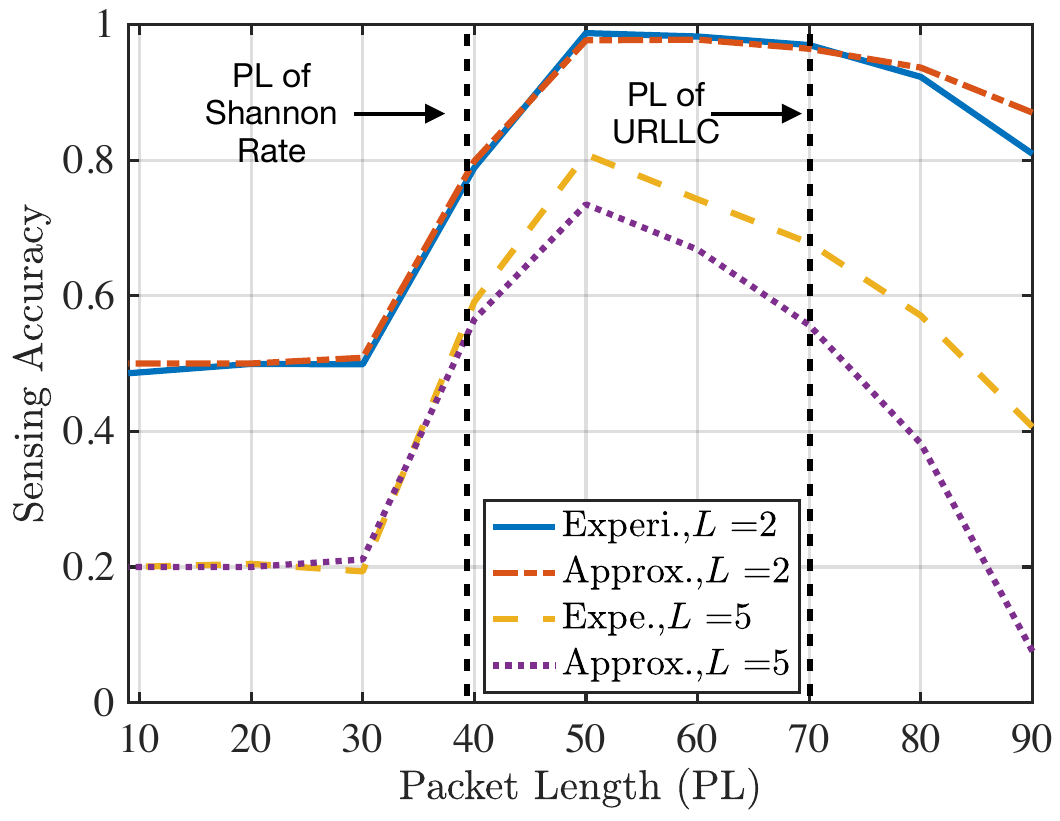}\label{Fig:MS_PL_ACC}}
\subfigure[Multi-view sensing (SNR = 10dB)]{
\includegraphics[width=0.45\columnwidth]{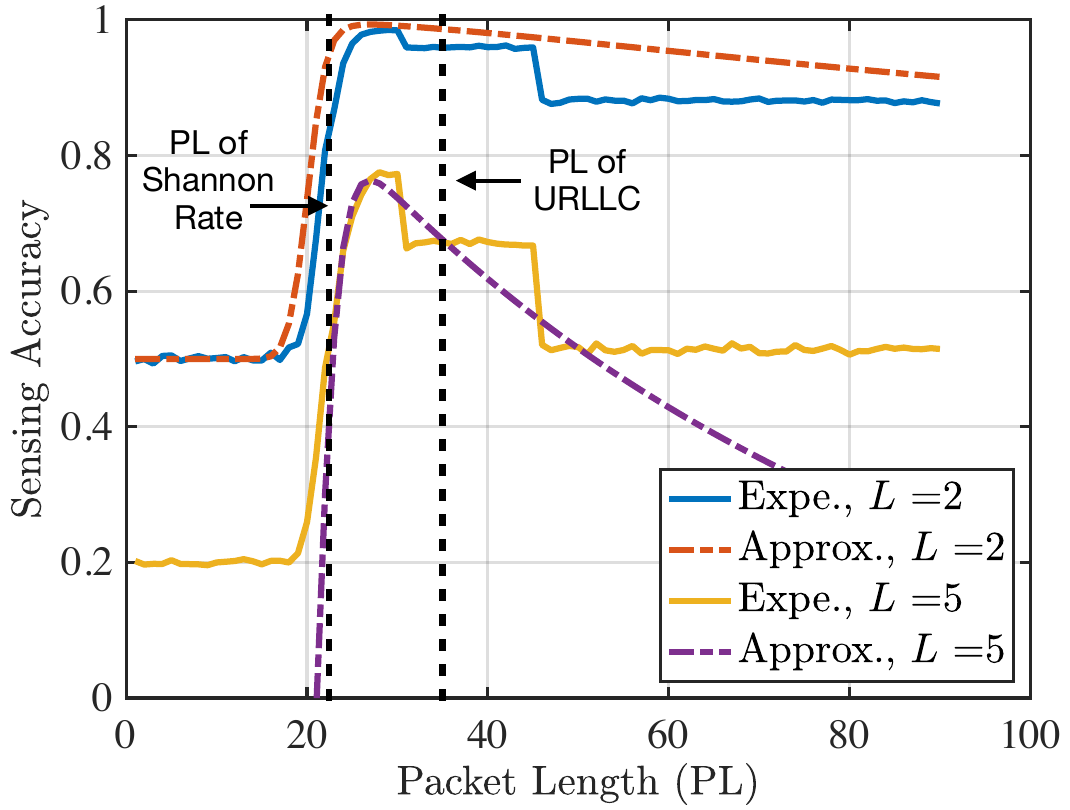}\label{Fig:MV_PL_ACC}}
\caption{Packet length versus sensing accuracy for linear classification.}
\label{Fig:SNR_PL_ACC}\vspace{-3mm}
\end{figure}

We start by presenting the derived reliability-view tradeoffs.
Fig. \ref{Fig:SNR_PL_ACC} compares the approximation of sensing accuracy in \eqref{sdp_discrete_rex} and \eqref{E2E_ACC_MV_approx}
with their ground truth in linear classification for a variable packet length.
For binary classification (i.e., $L=2$), the approximation can be seen to closely approximate the experimental performance across all packet lengths as also suggested by Remark \ref{Rem:Binary_Accuracy}.
For the multi-class classifier (i.e., $L=5$), the proposed approximation captures the relationship between the sensing accuracy and the packet length, suggesting that the approximation serves as a good surrogate for the sensing accuracy in the packet-length optimization.
Furthermore, the sensing accuracy is seen to increase and then decrease as the packet length increases, revealing the fundamental tradeoff between communication reliability and the number of views. 
Specifically, when the packet is long, the communication reliability is high, but the number of views is small (i.e., limiting the sensing quality). On the other hand, a short packet has low communication reliability but allows for many views (high sensing quality). Note that the packet lengths obtained using both URLLC and Shannon rate cannot attain the optimal performance. This is because the Shannon rate ignores the communication error caused by the finite blocklength, while URLLC ignores the sensing quality.

\begin{figure}[t!]
\centering
\subfigure[Multi-snapshot sensing ($L=5$)]{
\includegraphics[width=0.45\columnwidth]{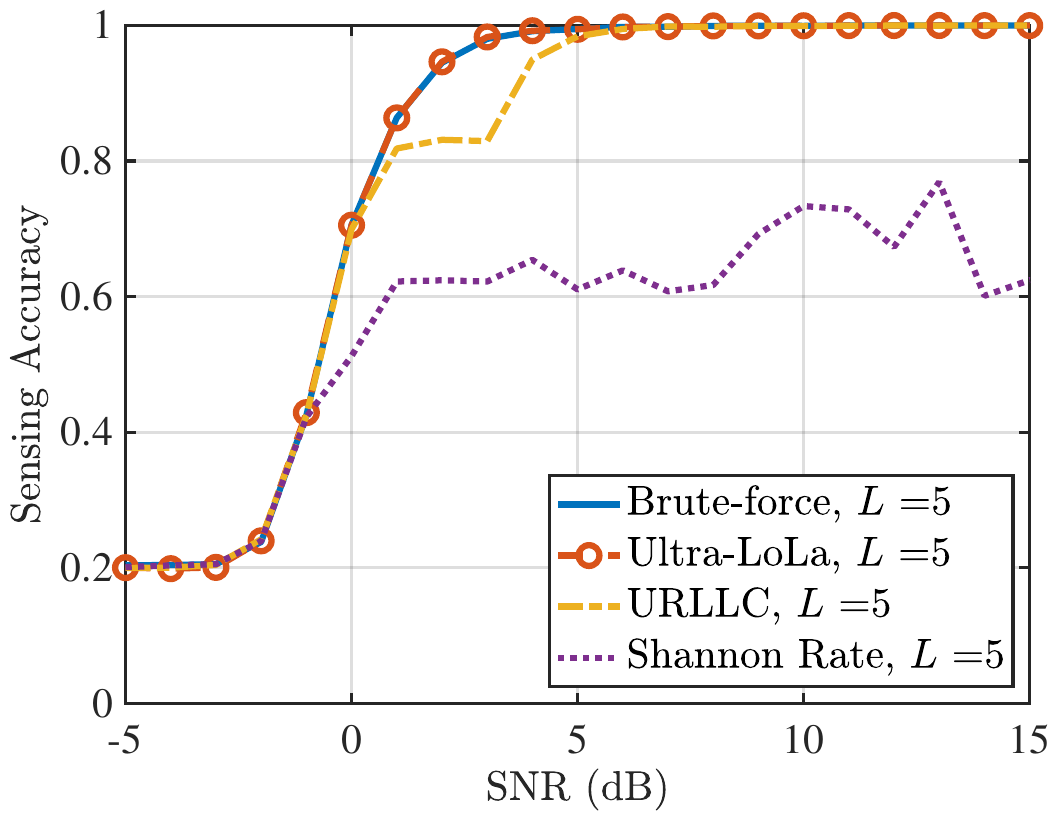}\label{Fig:MS_SNR_ACC_L5}}
\subfigure[Multi-snapshot sensing ($L=10$)]{
\includegraphics[width=0.45\columnwidth]{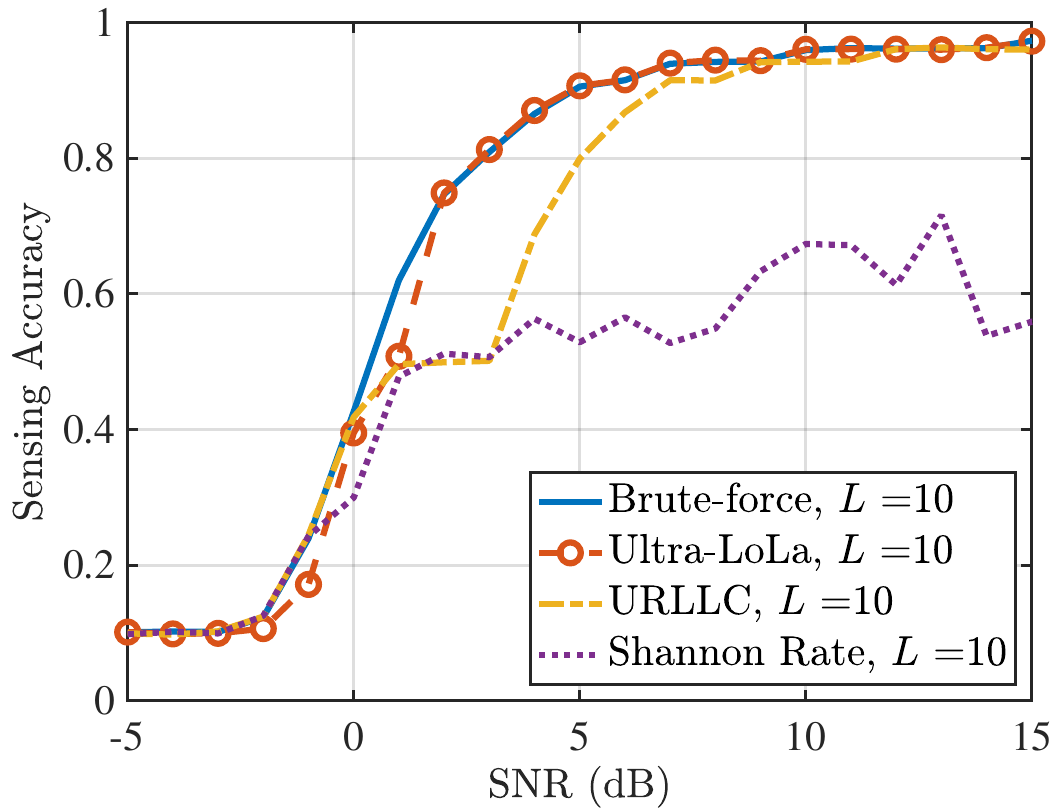}\label{Fig:MS_SNR_ACC_L10}}
\subfigure[Multi-view sensing ($L=5$)]{
\includegraphics[width=0.45\columnwidth]{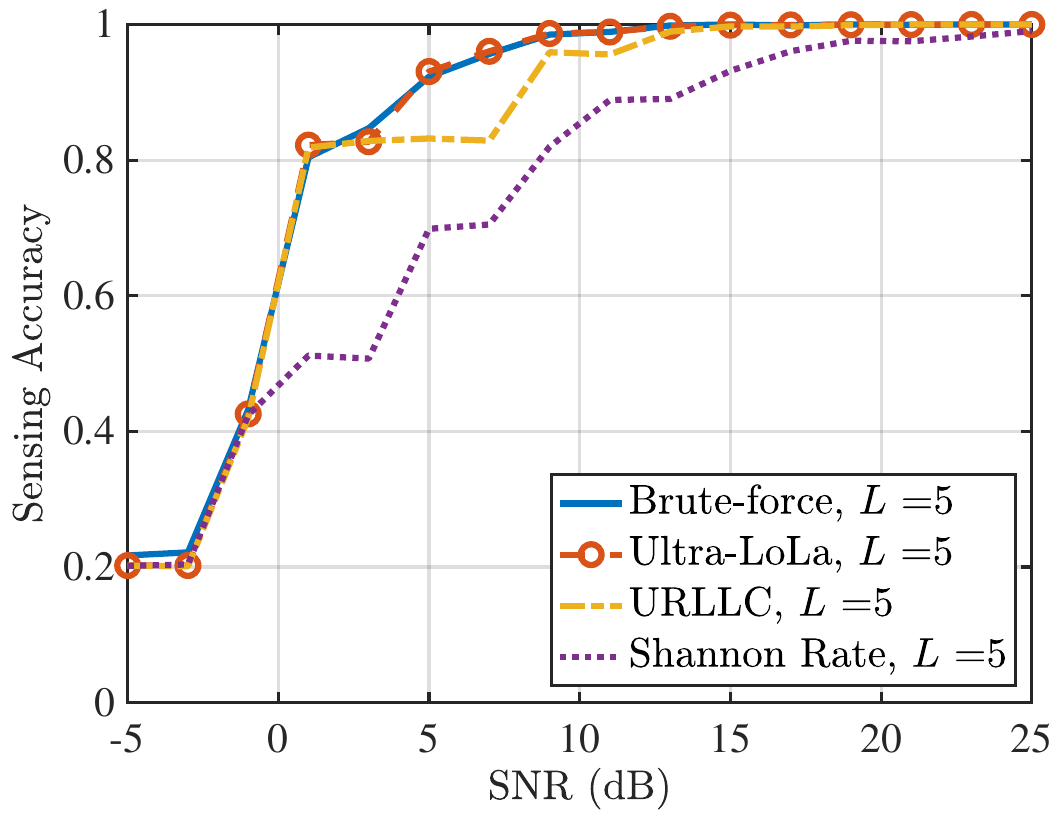}\label{Fig:MV_SNR_ACC_L5}}
\subfigure[Multi-view sensing ($L=10$)]{
\includegraphics[width=0.45\columnwidth]{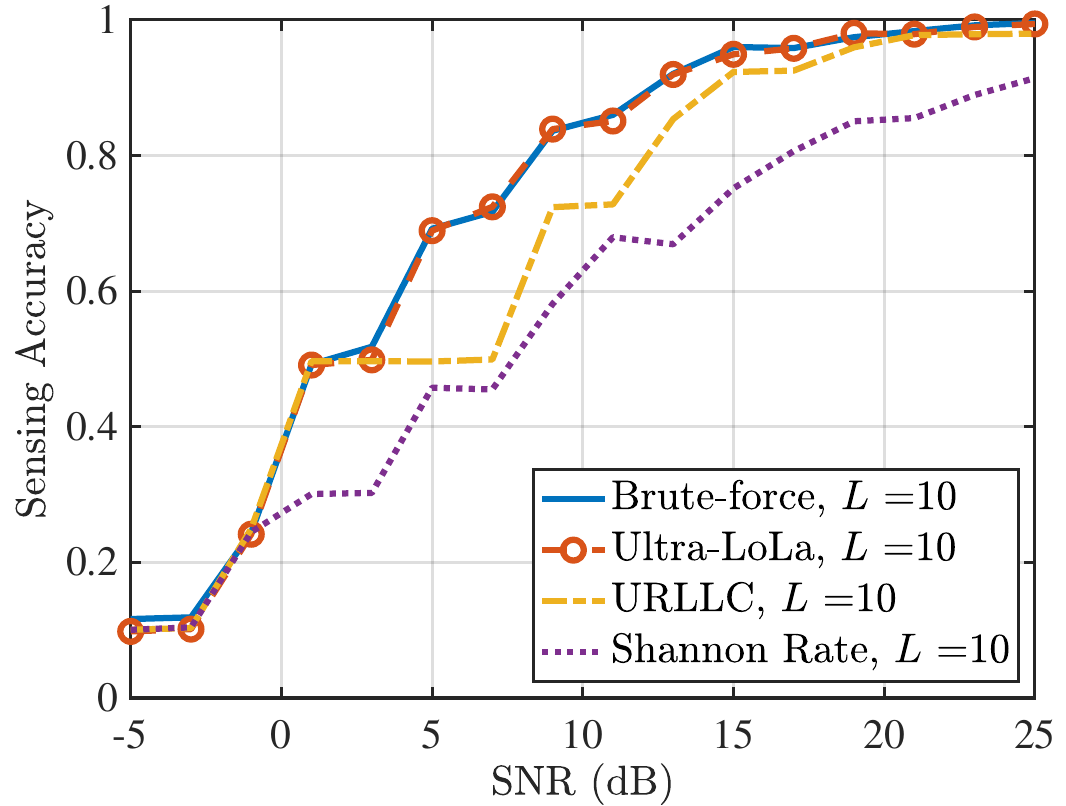}\label{Fig:MV_SNR_ACC_L10}}
\caption{Sensing accuracy comparison with benchmarks in the case of linear classification.}
\label{Fig:SNR_SNR_ACC}\vspace{-3mm}
\end{figure}

Fig. \ref{Fig:SNR_SNR_ACC} shows the sensing accuracy of ultra-LoLa inference as a function of SNR. 
The sensing accuracy of the ultra-LoLa packet length approaches the brute-force results, rising with the SNR and decreasing with the number of classes.
% This trend supports the theoretical analysis presented in Lemma \ref{Optimal_BL_Single_sensor} and  \ref{E2E_Acc_UP_MS}, indicating that E2E performance is a monotone increasing function of SNR and a decreasing function of the number of classes. 
The reasons are twofold:
First, a higher SNR leads to a lower decoding error rate, allowing the packet length to be reduced to include for more snapshots (or views) and the sensing quality.
Second, with a fixed feature size and covariance matrix, a larger number of classes gives a smaller pair-wise discrimination gain, which compromises classification accuracy.

In Fig. \ref{Fig:SNR_SNR_ACC}, we then compare ultra-LoLa inference with benchmarks in the case of linear classification.
The ultra-LoLa scheme achieves similar accuracies as benchmarks when SNR is lower than 0 dB  due to the high packet error probability. 
Compared with URLLC, the multi-view gain can be seen in the SNR range of 0-10 dB for multi-snapshot sensing and 3-20 dB for multi-view sensing. 
In this regime, ultra-LoLa inference includes more views at the cost of a higher decoding error rate,  while URLLC requiring ultra-reliable communication transmits fewer views due to a lower rate. 
For a high SNR, the required packet length of URLLC decreases, allowing more views to be included, which improves the sensing accuracy and results in a similar performance as ultra-LoLa inference.
Moreover, in the multi-snapshot case,
there is a significant performance gap between ultra-LoLa and
Shannon-rate schemes at high SNR. The reason is that the high
decoding error rate associated with the Shannon rate (around 0.5) becomes the bottleneck when all sensing information
is transmitted in a single packet. This phenomenon is less
significant in multi-view sensing due to the channel diversity obtained by having independently transmitted views \cite{jurdi2018outage}.

\begin{figure}[t!]
\centering
\subfigure[Multi-snapshot sensing]{
\includegraphics[width=0.45\columnwidth]{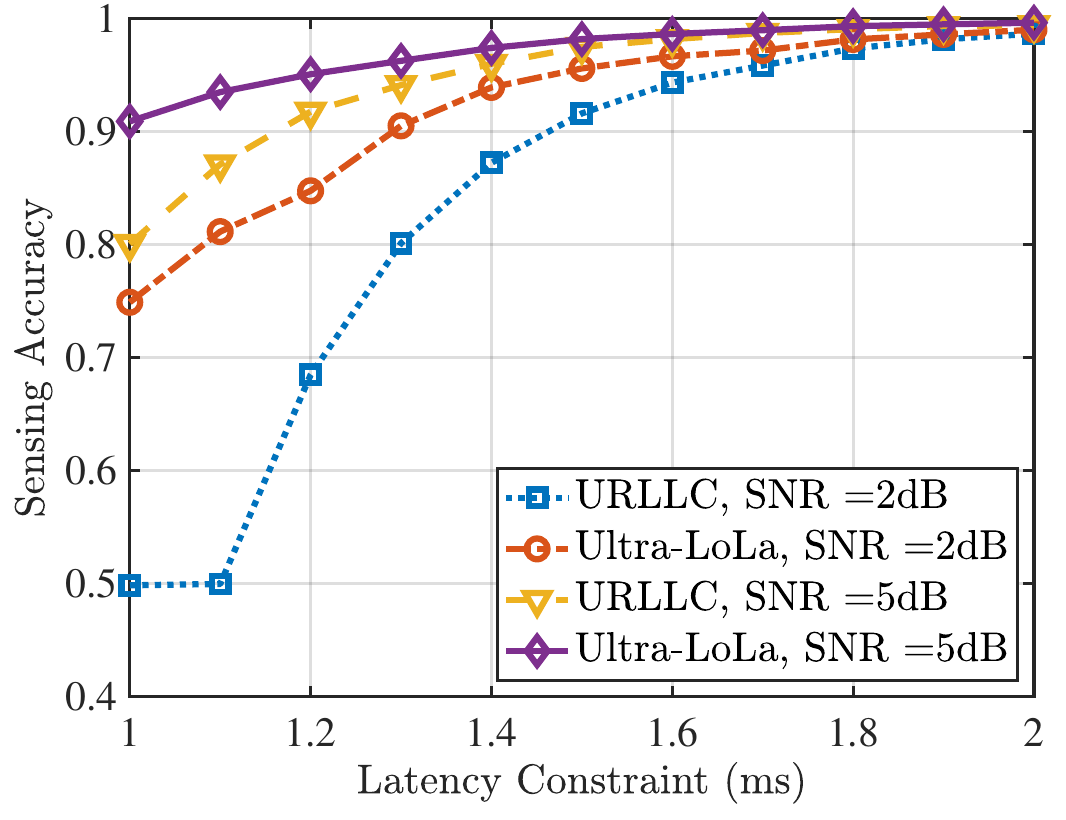}\label{Fig:LR_MS_LAT_ACC}}
\subfigure[Multi-view sensing]{
\includegraphics[width=0.45\columnwidth]{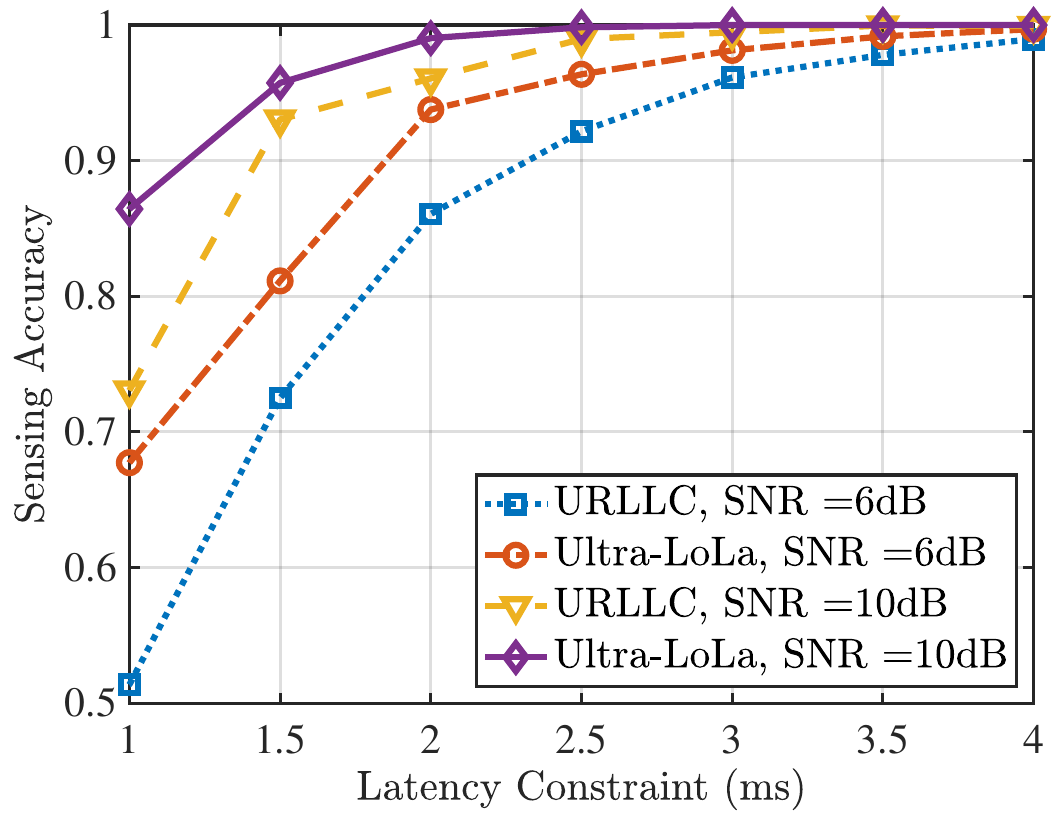}\label{Fig:LR_MV_LAT_ACC}}
\caption{Sensing accuracy versus latency constraints in the case of  10-class linear classifier.}
\label{Fig:LR_LAT_ACC}\vspace{-3mm}
\end{figure}

In Fig. \ref{Fig:LR_LAT_ACC}, 
the E2E performance of ultra-LoLa inference is measured by the E2E latency constraint and compared with the benchmarks.
Relaxing the latency constraint leads to an increase in sensing accuracy by allowing for longer packets and more views. 
It is observed that the gap between ultra-LoLa inference and the benchmarks is significant under strict latency constraints, showcasing the importance of balancing communication reliability and sensing quality. 
Additionally, for a given target sensing accuracy, ultra-LoLa inference attains  lower latency compared to other schemes.

\subsection{Ultra-LoLa Inference with MVCNN-Based Classification}
We now apply the proposed framework for MVCNN-based classification of a real dataset. We first optimize the packet length exploiting the insights from the linear classifier, and then compare the results to the benchmarks.

\subsubsection{Packet Length Optimization on Real Dataset}
 
To optimize the packet length of MVCNN, we employ the solution in Proposition \ref{Unique_Solution_MS} for multi-view sensing, while there is no theoretical solution for multi-snapshot sensing due to the unknown  minimum discrimination gain of the MVCNN. To tackle this challenge, we propose a lookup-table-based packet length optimization that exploits the insights from the  linear classifier.
Specifically, leveraging the linear relation between number of views and the packet length in \eqref{eq:num_observations}, we generate a table of packet lengths and corresponding training accuracy, denoted as $\Psi_{\sf{cnn}}(D)$, during training of the VGG16 using the training dataset.
The successful transmission probability can be expressed as a function of packet length with respect to the SNR, feature size, and quantization resolution using \eqref{succ_prob}.
Then we estimate the sensing accuracy of the MVCNN with the expectation of $\Psi_{\sf{cnn}}(D)$ over the decoding process, i.e., $\rho\Psi_{\sf{cnn}}(D)+\frac{(1-\rho)}{L}$.
In such manner, the packet length of ultra-LoLa inference can be estimated offline by an exhaustive search over the estimated sensing accuracy.

\begin{figure}
\centering
\subfigure[Multi-snapshot sensing]{
\includegraphics[width=0.45\columnwidth]{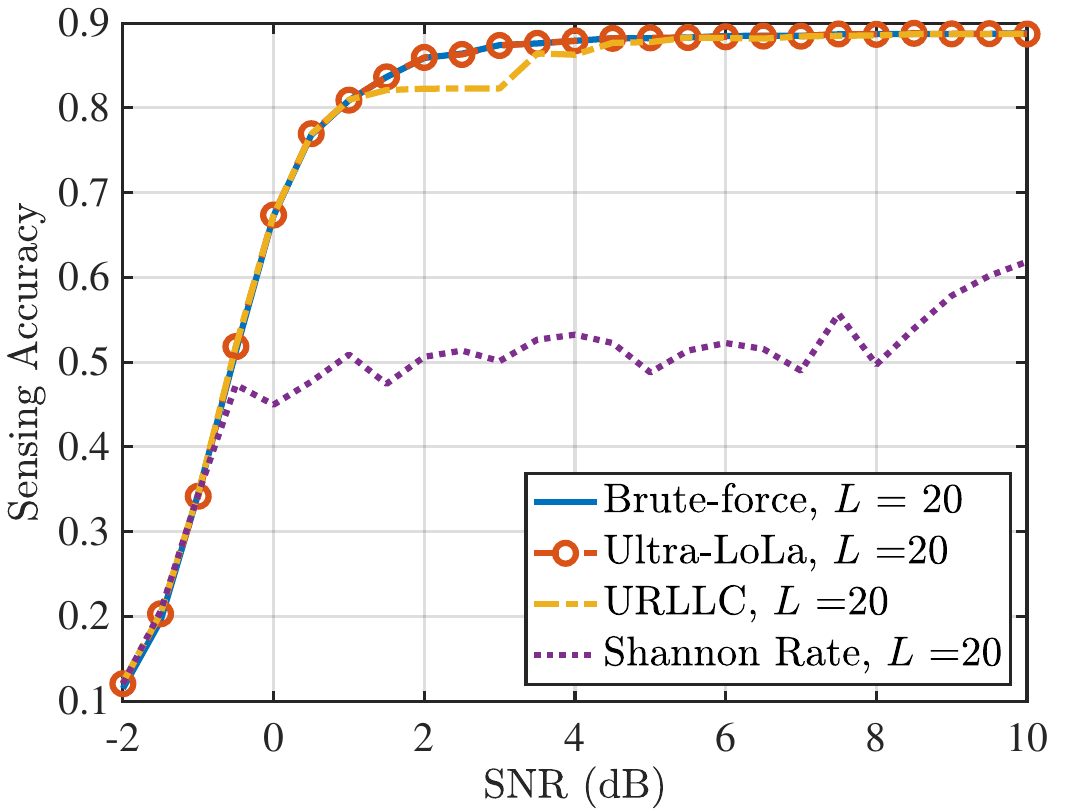}}
\subfigure[Multi-view sensing]{
\includegraphics[width=0.45\columnwidth]{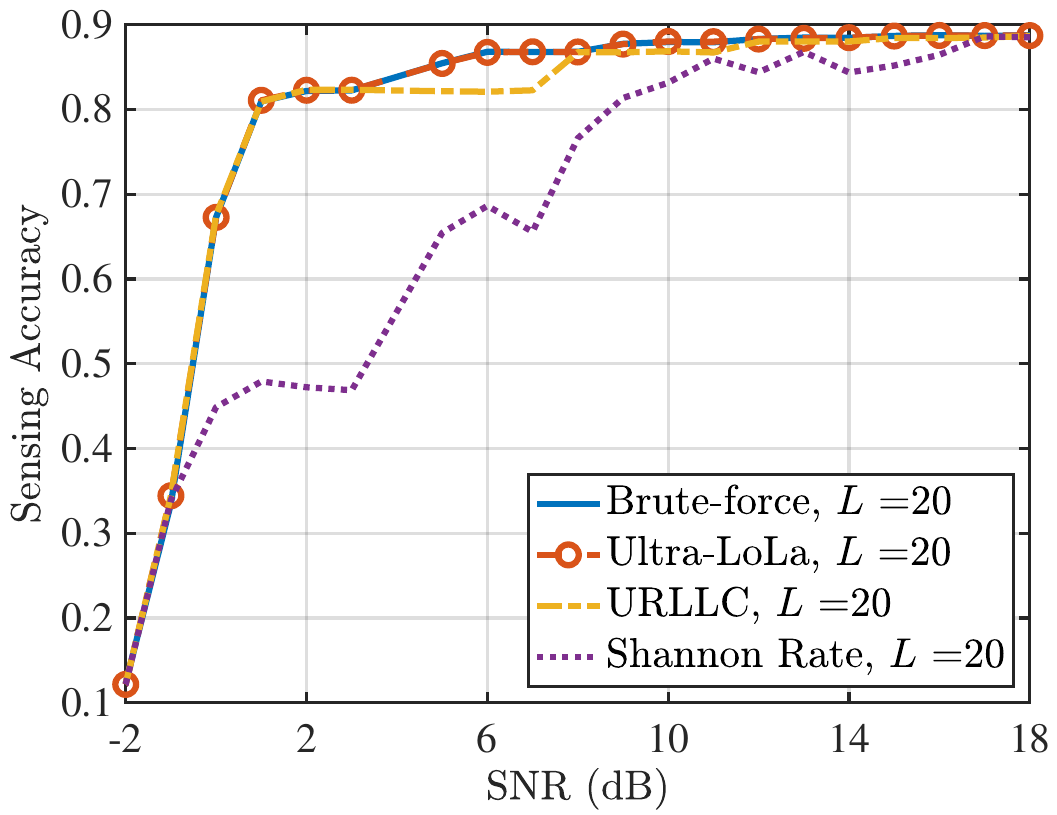}}
\caption{Sensing accuracy comparison with benchmarks in case of MVCNN-based classification.}
\label{Fig:CNN_SNR_ACC}
\vspace{-5mm}
\end{figure}

\subsubsection{Comparison with Benchmarks}
\begin{figure}
\centering
\subfigure[Multi-snapshot sensing]{
\includegraphics[width=0.45\columnwidth]{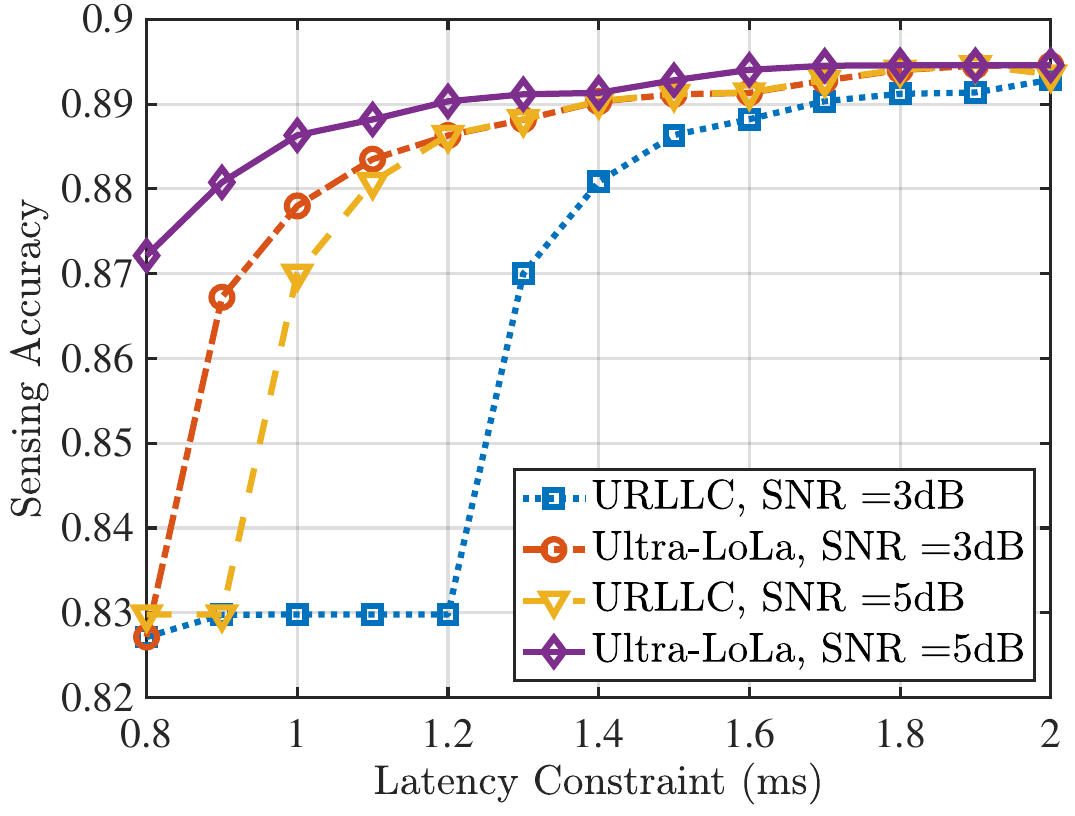}\label{Fig:CNN_MS_LAT_ACC}}
\subfigure[Multi-view sensing ]{
\includegraphics[width=0.45\columnwidth]{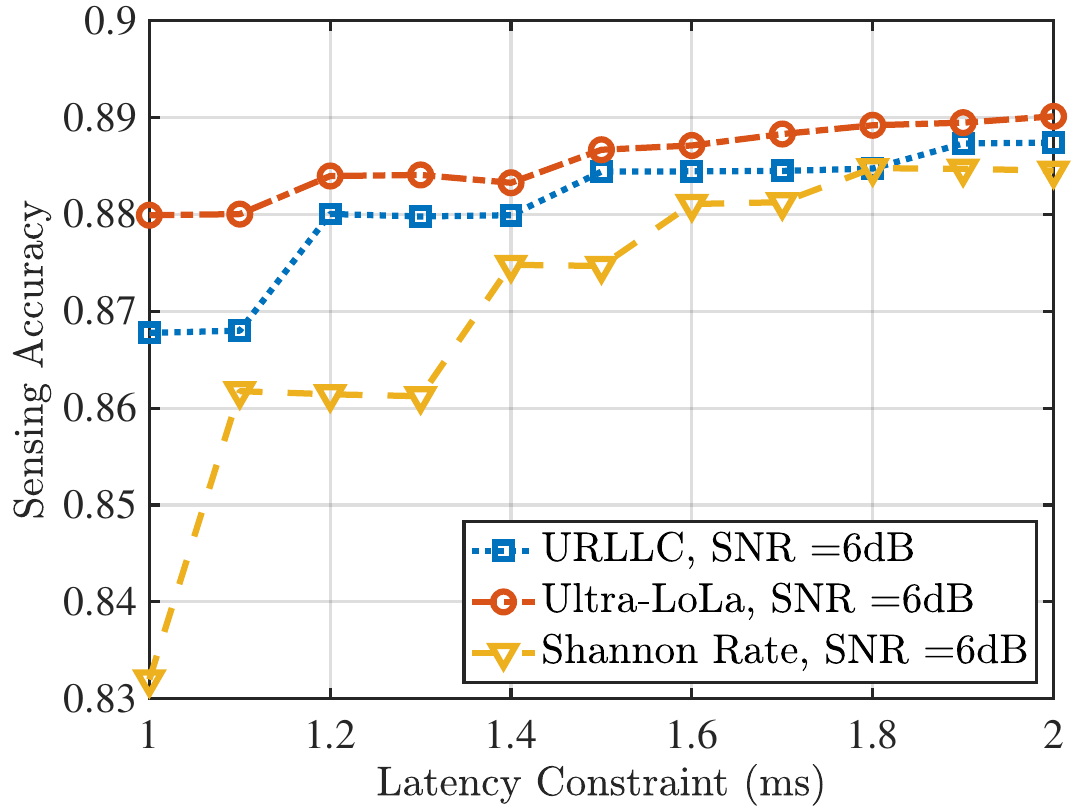}\label{Fig:CNN_MV_LAT_ACC }}
\caption{Sensing accuracy versus latency constraints in the case of  20-class MVCNN classifier.}
\label{Fig:CNN_LAT_ACC}
\vspace{-5mm}
\end{figure}

Figs. \ref{Fig:CNN_SNR_ACC} and \ref{Fig:CNN_LAT_ACC} demonstrate the performance of ultra-LoLa inference in comparison to the benchmarks in terms of sensing accuracy, latency and SNR.
As discussed for linear classification, similar trends can be observed in the MVCNN classifier.
For example, a multi-view gain between ultra-LoLa inference and URLLC can be observed at the SNR range of 1--4dB and 3--8dB for the two cases, respectively.
Among others, ultra-LoLa inference matches the brute-force results for all settings due to accurate prediction on the training dataset during model training.
Without this advantage, URLLC and Shannon rate fail to obtain the maximum performance at the same SNR settings.
In addition, as shown in Fig. \ref{Fig:CNN_LAT_ACC}, relaxing the latency constraint allows a larger packet length and a larger number of views, boosting the sensing accuracy as for linear classification.

\section{Conclusion}
\label{Conclusion}
In this paper, we have presented the framework of ultra-LoLa edge inference for latency-constrained distributed sensing.
The framework harnesses the interplay between SPT and accuracy improvements from multi-view sensing to meet a stringent deadline while boosting the E2E sensing performance.
Under the latency constraint, a fundamental tradeoff between communication reliability and the number of views, controlled by the packet length, is revealed and optimized.
The optimization is tackled by deriving accurate surrogate functions of the expressions for the E2E sensing accuracy.

This work is the first to study SPT in edge inference for distributed sensing, opening a number of opportunities for follow-up studies.
One direction is accuracy enhancement and latency reduction via joint optimization of packet length and coding rate.
This requires further work on characterizing the effects of feature dimensionality and quantization on packet length selection.
In terms of energy-efficient communication, it is interesting to design power control techniques that adapt to both fading channels and feature importance.
Additionally, the development of the ultra-LoLa inference to incorporate  AirComp warrants further investigation.

\vspace{-3mm}
\appendices

\section{Proof of Lemma \ref{inference_accuracy_LB}}
\label{Closed_form_IA}

Conditioned on successful transmission, the inference accuracy with $K$ observations is given as
\begin{equation}
\small
    A =\frac{1}{L}\sum_{\ell=1}^L\Pr(\hat{\ell}=\ell|\ell,K),
\end{equation}
where $\Pr(\hat{\ell}=\ell|\ell,K)$ is the correct classification probability given $K$ observations and ground-truth label $\ell$. %, i.e.,  $\mathbf{x}|\ell \sim \mathcal{N}(\boldsymbol{\mu}_\ell,C)$.
Following the observation model in \eqref{data_dis}, the conditional distribution of the fused feature vector $\overline{\mathbf{x}}$ defined in \eqref{feature_avepool} given $\ell$ is $\overline{\mathbf{x}}\sim\mathcal{N}(\boldsymbol{\mu}_{\ell},K^{-1}\mathbf{C})$.

Considering the Mahalanobis distance minimization classifier in \eqref{Maha_min_classifier}, the probability of correct classification is
\begin{comment}
% {\red
% I think it should be like this:
% \begin{equation}
% \label{Pr-1}
% \begin{split}
%      \Pr(\hat{\ell}=\ell|\ell) &= \Pr\left(\cap_{\ell'=1,\,\ell'\neq \ell}^L \left\{z_{\ell}(\mathbf{x}|\ell)<z_{\ell^\prime}(\mathbf{x}|\ell)\right\}     \right)\\
%      &= \Pr\left(\cap_{\ell'=1,\,\ell'\neq \ell}^L \left\{\delta_{\ell,\ell'}<0\right\}     \right)\\
%      &= 1-\Pr\left(\cup_{\ell'=1,\,\ell'\neq \ell}^L \left\{\delta_{\ell,\ell'}<0\right\}     \right)\\
%      &\ge 1-\sum_{\ell'=1,\,\ell'\neq \ell}^L\Pr\left(\delta_{\ell,\ell'}<0\right),
% \end{split}
% \end{equation}
% Alternative simple bound (Fréchet inequalities):
% \begin{align}
%     \Pr\left(\cap_{\ell'=1,\,\ell'\neq \ell}^L \left\{\delta_{\ell,\ell'}<0\right\}     \right)&\ge \sum_{\ell'=1,\,\ell'\neq \ell}^L\Pr\left(\delta_{\ell,\ell'}<0\right)-(L-2)\\
%     &\ge (L-1)\min_{\ell'=1,\ldots,L;\,\ell'\neq \ell}\Pr\left(\delta_{\ell,\ell'}<0\right)-(L-2)
% \end{align}
% Another simple one is to bound the union bound sum:
% \begin{align}
%     \Pr\left(\cap_{\ell'=1,\,\ell'\neq \ell}^L \left\{\delta_{\ell,\ell'}<0\right\}     \right)&\ge 1-(L-1)\max_{\ell'=1,\ldots,L;\,\ell'\neq \ell}\Pr\left(\delta_{\ell,\ell'}<0\right)
% \end{align}
% }
\end{comment}
\begin{equation}
    \Pr(\hat{\ell}=\ell|\ell,K) = \Pr\left(\cap_{\ell'=1,\,\ell'\neq \ell}^L \left\{z_{\ell}(\overline{\mathbf{x}}|\ell)<z_{\ell^\prime}(\overline{\mathbf{x}}|\ell)\right\}     \right),
\end{equation}
where $z_{\ell'}(\overline{\mathbf{x}}|\ell)=(\overline{\mathbf{x}}-\boldsymbol{\mu}_{\ell'})^{\sf{T}}\mathbf{C}^{-1}(\overline{\mathbf{x}}-\boldsymbol{\mu}_{\ell'})$ is the squared Mahalanobis distance between $\overline{\mathbf{x}}$ and class $\ell'$ given that the ground-truth class is $\ell$.
Defining the differential squared Mahalanobis distance as $\delta_{\ell,\ell^\prime}=z_{\ell}(\overline{\mathbf{x}}|\ell)-z_{\ell^\prime}(\overline{\mathbf{x}}|\ell)$, we have
\begin{equation}
\label{Pr-lb}
\begin{split}   \small
\Pr(\hat{\ell}=\ell|\ell,K)
 &= \Pr\left(\cap_{\ell'=1,\,\ell'\neq \ell}^L \left\{\delta_{\ell,\ell'}<0\right\}     \right)\\
 &= 1-\Pr\left(\cup_{\ell'=1,\,\ell'\neq \ell}^L \left\{\delta_{\ell,\ell'}\geq0\right\}     \right)\\
 &\ge 1-\sum_{\ell'=1,\,\ell'\neq \ell}^L\Pr\left(\delta_{\ell,\ell'}\geq0\right)\\
  &\geq 1- (L-1) \max_{\forall \ell,\ell'}\Pr\left(\delta_{\ell,\ell'}\geq0\right).
\end{split}
\end{equation}
Denoting by $\overline{\mathbf{x}}_{\ell}(n)$ the $n$-th entry of $\overline{\mathbf{x}}$ when the ground-truth class is $\ell$, the differential squared Mahalanobis distance can be further expressed as $\delta_{\ell,\ell^\prime}=z_{\ell}(\overline{\mathbf{x}}|\ell)-z_{\ell^\prime}(\overline{\mathbf{x}}|\ell) = \sum_{n=1}^N \delta_{\ell,\ell^\prime}(n)$, 
% \begin{equation}
%    \begin{split} \small
% &\delta_{\ell,\ell^\prime}=z_{\ell}(\overline{\mathbf{x}}|\ell)-z_{\ell^\prime}(\overline{\mathbf{x}}|\ell)
%         % =&\sum_{n=1}^N \frac{(\overline{\mathbf{x}}_\ell(n)-\boldsymbol{\mu}_{\ell}(n))^2}{C_{n,n}}-\frac{(\overline{\mathbf{x}}_\ell(n)-\boldsymbol{\mu}_{\ell^\prime}(n))^2}{C_{n,n}}\\
% %         =&\sum_{n=1}^N \underbrace{\frac{2}{C_{n,n}}(\boldsymbol{\mu}_{\ell^\prime}(n)-\boldsymbol{\mu}_{\ell}(n))
% % \left(\overline{\mathbf{x}}_\ell(n)-\frac{1}{2}\left(\boldsymbol{\mu}_{\ell}(n)+\boldsymbol{\mu}_{\ell^\prime}(n)\right)\right)}_{\delta_{\ell,\ell^\prime}(n)}\\
%         = \sum_{n=1}^N \delta_{\ell,\ell^\prime}(n),
%    \end{split}
% \end{equation}
where $\delta_{\ell,\ell^\prime}(n)$ is the differential distance of the $n$-th dimension, given as
\begin{equation}\small
    \delta_{\ell,\ell^\prime}(n)=\frac{2}{C_{n,n}}(\boldsymbol{\mu}_{\ell^\prime}(n) - \boldsymbol{\mu}_{\ell}(n))
\left(\overline{\mathbf{x}}_\ell(n)-\frac{1}{2}\left(\boldsymbol{\mu}_{\ell}(n)+\boldsymbol{\mu}_{\ell^\prime}(n)\right)\right).
\end{equation}
Note that $\delta_{\ell,\ell^\prime}(n)\sim \mathcal{N}(-g_{\ell,\ell^\prime}(n),4K^{-1}g_{\ell,\ell^\prime}(n))$ is Gaussian consdering the linear mapping from
$\overline{\mathbf{x}}_\ell(n)$.
%Since $\overline{\mathbf{x}}_\ell(n)\sim \mathcal{N}(\boldsymbol{\mu}_\ell(n),K^{-1}C_{n,n})$,  the differential distance $\delta_{\ell,\ell^\prime} (n)$ is also Gaussian.
 %\begin{equation}
 %    \delta_{\ell,\ell^\prime} (n) \sim \mathcal{N}(-g_{\ell,\ell^\prime}(n),4K^{-1}g_{\ell,\ell^\prime}(n)),
 %\end{equation}
 %where $g_{\ell,\ell^\prime}(n)=\frac{(\boldsymbol{\mu}_{\ell^\prime}(n)-\boldsymbol{\mu}_{\ell}(n))^2}{C_{n,n}}$ is the $n$-th dimension discrimination gain between class $\ell$ and class $\ell^\prime$.
 % \begin{equation}
 %     g_{\ell,\ell^\prime}(n)=\frac{(\boldsymbol{\mu}_{\ell^\prime}(n)-\boldsymbol{\mu}_{\ell}(n))^2}{C_{n,n}}.
 % \end{equation}
Since the covariance matrix, $\mathbf{C}$, of the observations is diagonal, $\delta_{\ell,\ell^\prime}$ is a sum of $N$ independent Gaussian random variables, distributed as
 \begin{equation}
     \delta_{\ell,\ell^\prime} \sim \mathcal{N}(-g_{\ell,\ell^\prime},4K^{-1}g_{\ell,\ell^\prime}),
 \end{equation}
 where $g_{\ell,\ell^\prime}=\sum_{n=1}^N g_{\ell,\ell^\prime}(n)=[\boldsymbol{\mu}_\ell-\boldsymbol{\mu}_{\ell^\prime}]^{\sf{T}}\mathbf{C}^{-1}[\boldsymbol{\mu}_\ell-\boldsymbol{\mu}_{\ell^\prime}]$ is the total discrimination gain between the two classes.

The probability that distance from ground truth label $\ell$ is greater than that of $\ell^\prime$, i.e, $\Pr(\delta_{\ell,\ell^\prime}\ge0)$, can then be given as
 \begin{equation}
 \begin{split}
 \label{Pr-2}
    \small
\Pr(\delta_{\ell,\ell^\prime}\ge 0) &=  \int_{0}^\infty \mathcal{N}(\delta_{\ell,\ell^\prime};-g_{\ell,\ell^\prime},4K^{-1}g_{\ell,\ell^\prime}) d\delta_{\ell,\ell^\prime} \\
      &=Q\left(\frac{\sqrt{Kg_{\ell,\ell^\prime}}}{2}\right).
  \end{split}
 \end{equation}
Substituting \eqref{Pr-2} into \eqref{Pr-lb}, we obtain the lower bound of inference accuracy, given as
\begin{equation}
    \begin{split}
 A =\frac{1}{L}\sum_{\ell=1}^L \Pr(\hat{\ell}=\ell|\ell)
  \geq 1-(L-1)Q\left( \frac{\sqrt{Kg_{\min}}}{2}\right),
\end{split}
\end{equation}
where $g_{\min}=\min\{ g_{\ell_1,\ell_2}|\ell_1,\ell_2\in \{1,2,\dots,L\}, \ell_1\neq\ell_2 \}$ is the minimum discrimination gain between any two classes. The proof is completed by noting that the accuracy cannot be smaller than $1/L$ (achieved by random guessing).
This completes the proof.

\section{Proof of Proposition  \ref{Optimal_BL_Single_sensor}}
\label{Proof_Existence_Optimal_BL}

We first show that $\tilde{\nu}_{\sf{ms}}(D)$ is a log-concave function over the continuous (relaxed) domain $D \in [1,D_{\max}]$, and then establish the optimality conditions for the original problem with $D\in \{1,\ldots,D_{\max}\}$.

Define $f_I(D)=\sigma(\psi_I(D))$ and $f_T(D)=\sigma(\psi_T(D))$, which are positive functions of $D \in [1,D_{\max}]$. $\tilde{\nu}_{\sf{ms}}(D)$ can then be expressed as
\begin{equation}
     \tilde{\nu}_{\sf{ms}}(D)=\left(f_I(D)-\frac{L-1}{L}\right) f_T(D).
\end{equation}
The functions $\psi_I(x)$ and $\psi_T(x)$ are concave, which can be seen from their negative second-order derivatives. Specifically, for $D \in [1,D_{\max}]$
\begin{equation}
\begin{split}
\small
     \psi^{\prime\prime}_I(D)& = -\frac{G\eta}{4\Delta_S^2 B_W^2}\left(\frac{T}{\Delta_S}-\frac{D}{\Delta_S B_W } \right)^{-\frac{3}{2}}< 0,\\
     \psi^{\prime\prime}_T(D) & = -\frac{\ln(2)\eta}{4\sqrt{V}} D^{-\frac{5}{2}}\left( \log_2(1+\gamma)D+3Q_B N\right)<0.
\end{split} 
\end{equation}
Similarly, the sigmoid function $\sigma(x)=\frac{1}{1+\exp(-x)}$ is a log-concave and monotonically increasing function of $x$, and thus $f_I(D)$ and $f_T(D)$ can be proved to be log-concave functions of $D$~\cite{boyd2004convex}. 
Building on the constraint  that classification accuracy is larger than that of random guessing, i.e.,
\begin{equation}
\label{eq:D_max_con}
    Q\left(G\sqrt{\frac{T}{\Delta_S}-\frac{D}{B_W \Delta_S}}\right)\leq 0.
\end{equation}
The packet length is constrained to at most $D_{\max}$  given in \eqref{D_max_upper}.
Under this constraint, we have $f_I(D)-\frac{L-1}{L}>0, \forall D\in [1, D_{\max}] $, and thus
$f_I(D)-\frac{L-1}{L}$ is also a log-concave function in $ D\in [1, D_{\max}]$ (see, e.g., \cite[Exercise 3.48]{boyd2004convex}). 
The products of two log-concave functions is also a log-concave function, and thus $\tilde{\nu}_{\sf{ms}}(D)$ is a log-concave function of $D\in [1,D_{\max}]$. 
Then considering the monotonicity of $\ln(x)$, the optimal packet length $\Tilde{D}^*$ of problem \eqref{prob:max_E2E_acc_MC} can be given as
\begin{equation}
\begin{split}
     \Tilde{D}^*&=\argmax_D \tilde{\nu}_{\sf{ms}}(D)\\
     &=\argmax_D \ln(\tilde{\nu}_{\sf{ms}}(D))\\
     &=\{D| \tilde{\nu}'_{\sf{ms}}(D)=0,D\in  [1,D_{\max}] \} .
\end{split}
\end{equation}
Guaranteeing the maximum point reside in $[1,D_{\max}]$ is equivalent to the existence of the zero of continuous function $ \tilde{\nu}'_{\sf{ms}}(D), \forall D \in [1,D_{\max}]$. 
Thus, $\tilde{\nu}'_{\sf{ms}}(D)$ evaluated at the endpoints must satisfy
\begin{equation}
\label{ineq1}
{\Tilde{{\nu}}'_{\sf{ms}}(1)\cdot\Tilde{\nu}}'_{\sf{ms}}(D_{\max})<0, 
\end{equation}
where $\Tilde{{\nu}}'_{\sf{ms}}(D)$ is the first derivative given as 
\begin{equation}
\begin{split} \Tilde{{\nu}}'_{\sf{ms}}(D)=  f_{\sf{pos}}(D)f_{\bf{ms}}(D),
\end{split}
\end{equation}
with $f_{\sf{pos}}(D)= -\sigma(\psi_T(D)) (1-\sigma(\psi_T(D)))  \left(\sigma(\psi_I(D))-\frac{L-1}{L}  \right)\psi'_I(D) \geq 0, \forall D\geq 0$, and 
$f_{\bf{ms}}(D)$ being a function of $D$ determining the sign of $\Tilde{{\nu}}'_{\sf{ms}}(D)$, given as
\begin{equation}
\begin{split}
     f_{\bf{ms}}(D)&=-\frac{ 1+e^{\psi_T(D)} }{(1+e^{\psi_I(D)})(1/L-e^{-\psi_I(D)})}-\underbrace{\frac{\psi'_T(D)}{\psi'_I(D)}}_{(a)}.\\
\end{split}
\end{equation}
Note that $(a)$ can be expressed as
\begin{equation}
\begin{small}
    \frac{\psi'_T(D)}{\psi'_I(D)}=-C_{\bf{ms}} \sqrt{\frac{B_WT}{D}-1}\left( \frac{\log_2(1+\gamma)D+Q_BN}{D}\right),
    \end{small}
\end{equation}
where $C_{\bf{ms}}=\frac{\ln(2)}{G}\sqrt{\frac{\Delta_S B_W}{V}}$ is a positive constant.
When the inequality $  f_{\bf{ms}}(1)\cdot f_{\bf{ms}}(D_{\max} )<0$ holds, the optimal packet length $\Tilde{D}^*$ is obtained by finding the zero point of $f_{\bf{ms}}(\Tilde{D}^*)=0$.
Since the problem is log-convex, the optimal integer packet length, denoted as $D^*$, can be obtained by simply inspecting the nearest (feasible) integers less than or greater than $\Tilde{D}^*$ and picking the one that maximizes the objective function $\Tilde{\nu}_{\sf{ms}}(D)$.
Otherwise, if $  f_{\bf{ms}}(1)\cdot f_{\bf{ms}}(D_{\max} )\ge 0$, the optimal packet length is one of the endpoints, i.e., $D^*= \argmax_{D\in\{1,D_{\max}\}}\tilde{\nu}_{\sf{ms}}(D)$.
This completes the proof.

\vspace{-3mm}
\section{Proof of Proposition \ref{Unique_Solution_MS}}
\label{Proof_Unique_Solution_Prob}

First, to find the  maximum of $\Tilde{\nu}_{\sf{mv}}(D)$ we investigate its monotonicity through its derivative
\begin{equation}
\begin{split}
    \label{FD_nu(D)}    \Tilde{\nu}_{\sf{mv}}^\prime(D)=&C_{\sf{pos}} \underbrace{\left[ -(1+\exp(\psi_T(D)))+D\psi^\prime_T(D) \right] }_{g(D)},
\end{split}
\end{equation}
where $C_{\sf{pos}}=\frac{(T-\Delta_S)B_W\xi_a\sigma(\psi_T(D))} {D^2  (1+\exp{(\psi_T(D))})}>0$ is a positive constant.
The monotonicity of $\Tilde{\nu}_{\sf{mv}}(D)$ thus depends on the sign of $g(D)$ in \eqref{FD_nu(D)}. $g(D)$ can be expanded as
\begin{equation}
\begin{split}
\small
\label{g(D)}
    g(D)=&-\exp{\left[ \frac{\ln(2)\eta}{\sqrt{V}}\left(\frac{\log_2(1+\gamma)D -Q_BN}{\sqrt{D}}\right)\right]}\\
    &+\frac{\ln(2)\eta}{2\sqrt{V}}\left(\log_2(1+\gamma)\sqrt{D} +\frac{Q_BN}{\sqrt{D}}\right)-1.
    \end{split}
\end{equation}

Define the positive constant 
$\alpha=\frac{\ln(2)\eta}{\sqrt{V}}\log_2(1+\gamma)$ such that $g(D)=q\left(\alpha\sqrt{D} \right)$ where
\begin{equation}
\begin{split}
     q(x) = \frac{1}{2}\left(x+\frac{\omega}{x}\right)-\left(\exp\left(x-\frac{\omega}{x}\right)+1\right),
\end{split}
\end{equation}
for $x>0$ with $\omega=\frac{\alpha^2Q_BN}{\log_2(1+\gamma)}>1$ being a positive constant.
Then the first derivative of $g(D)$ can be obtained as
\begin{equation}\label{eq:gprime}
g'(D)=\left(\frac{\alpha}{2\sqrt{D}}\right)q'(\alpha \sqrt{D}),
\end{equation}
where $q'(x)$ is the first derivative of $q(x)$ with respect to $x$, given as
\begin{equation}
\begin{split}
      q^\prime(x)&=\frac{1}{2x^2}\left[x^2-\omega-2(x^2+\omega)\exp(x-\frac{\omega}{x})\right]\\
      & =\begin{cases}
          < 0, & x\in (0,\sqrt{\omega})\\
          < \frac{1}{2x^2}\left[-x^2-3\omega\right]< 0,& x\in [\sqrt{\omega},\infty).
      \end{cases}
\end{split}
\end{equation}
Since $q'(x)<0$ for $x>0$, this holds for $g'(x)$ as well.
Thus $g(D)$ is monotonically decreasing in $D\in [1,D_{\max}]$.

Since $g(D)$ is monotonically decreasing, the value $\Tilde{D}^*$ that maximizes $\Tilde{\nu}_{\sf{mv}}^\prime(D)$ either satisfies $g(\Tilde{D}^*)=0$, or is one on the boundary of the domain of $D$. The existence of a $D$ that satisfies $g(D)=0$ is guaranteed if and only if $g(1)g(D_{\max})< 0$, such that the boundary points of $g(D)$ have different signs. This condition is the same as the one in \eqref{condtion_optimalBL_MS}.
%Considering the continuity and monotonicity of $g(D)$, there must exist a unique value $\Tilde{D}^*$ satisfying $g(\Tilde{D}^*)=0$ if and only if $g(1)g(D_{\max})< 0$. such that $g(\Tilde{D}^*)=0, g(D)> 0 ,\forall D \in[1,\Tilde{D}^*) $ and $g(D)< 0 ,\forall D \in[\Tilde{D}^*,D_{\max})$ hold.
%Given the same sign of $\Tilde{\nu}'_{\sf{mv}}(D)$ and $g(D)$, $\Tilde{\nu}_{\sf{mv}}(D), D\in [1,D_{\max}] $ is a unimodal function with a maximum.
Thus, assuming that the condition in \eqref{condtion_optimalBL_MS} is satisfied, the optimal packet length, denoted as $D^*$, is unique and given by the solution to transcendental equation $q(\hat{x})=0$, i.e.,
\begin{equation}
        \hat{x}+\frac{\omega}{\hat{x}}-2\exp{\left( \hat{x}- \frac{\omega}{\hat{x}}\right)}=2.
\end{equation}
Finally, the effects of rounding on the choice of integer packet length is the same as discussed in Appendix \ref{Proof_Existence_Optimal_BL}, and thus is omitted.
This completes the proof.

\bibliography{Main_Accepted}
\bibliographystyle{IEEEtran}

\begin{IEEEbiography}[{\includegraphics[width=1in,height=1.25in, clip,keepaspectratio]{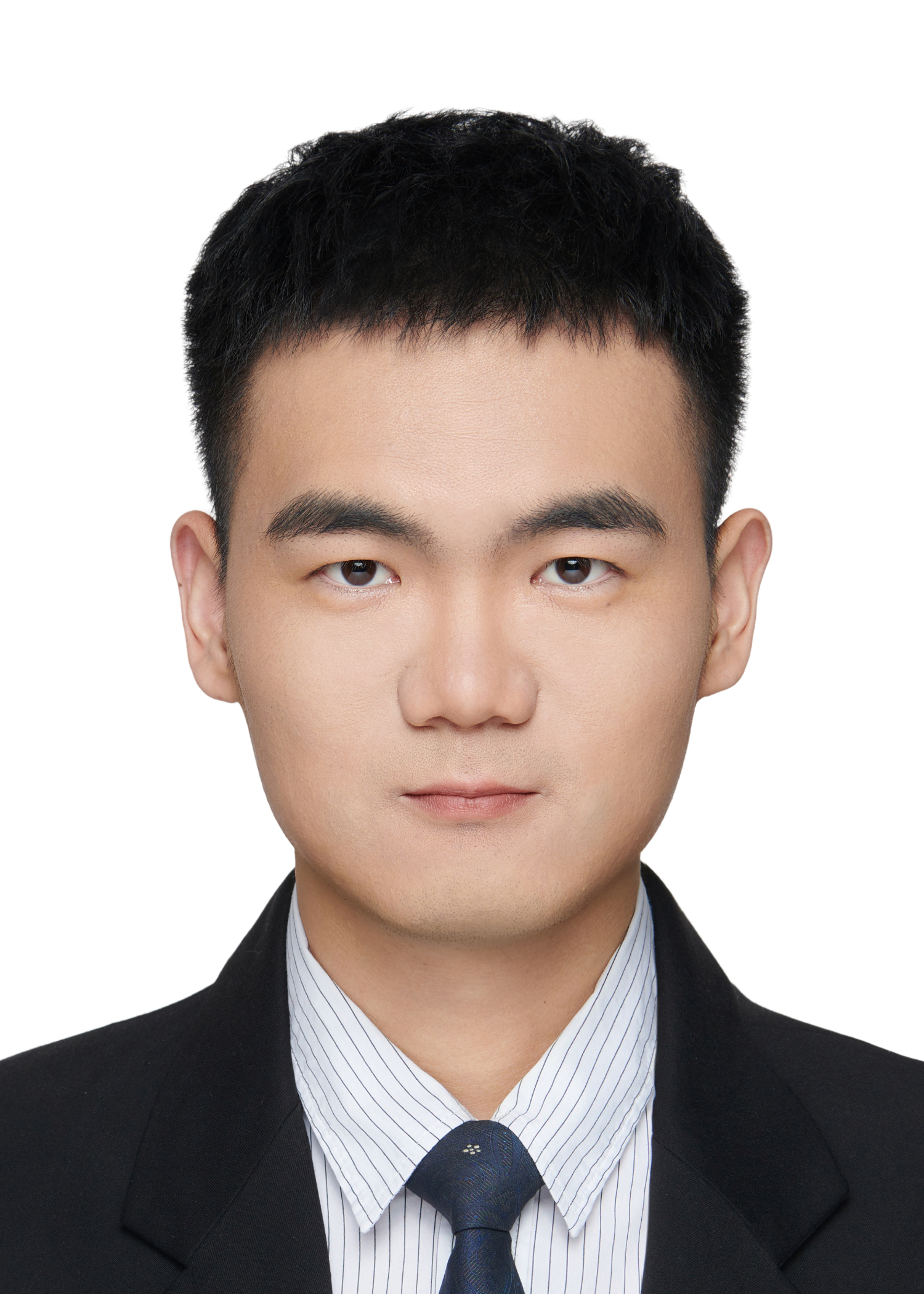}}]{Zhanwei Wang} (Graduate Student Member, IEEE) received the B.Eng. degree in Information Engineering and the M.Eng. degree in Information and Communication Engineering from Xidian University, Xi’an, China, in 2018 and 2021, respectively. He is currently pursuing a Ph.D. degree in the Department of Electrical and Electronic Engineering at the University of Hong Kong. His research interests include wireless communications, edge intelligence, distributed sensing, and atomic receiver.
\end{IEEEbiography}

\begin{IEEEbiography}[{\includegraphics[width=1in,height=1.25in, clip,keepaspectratio]{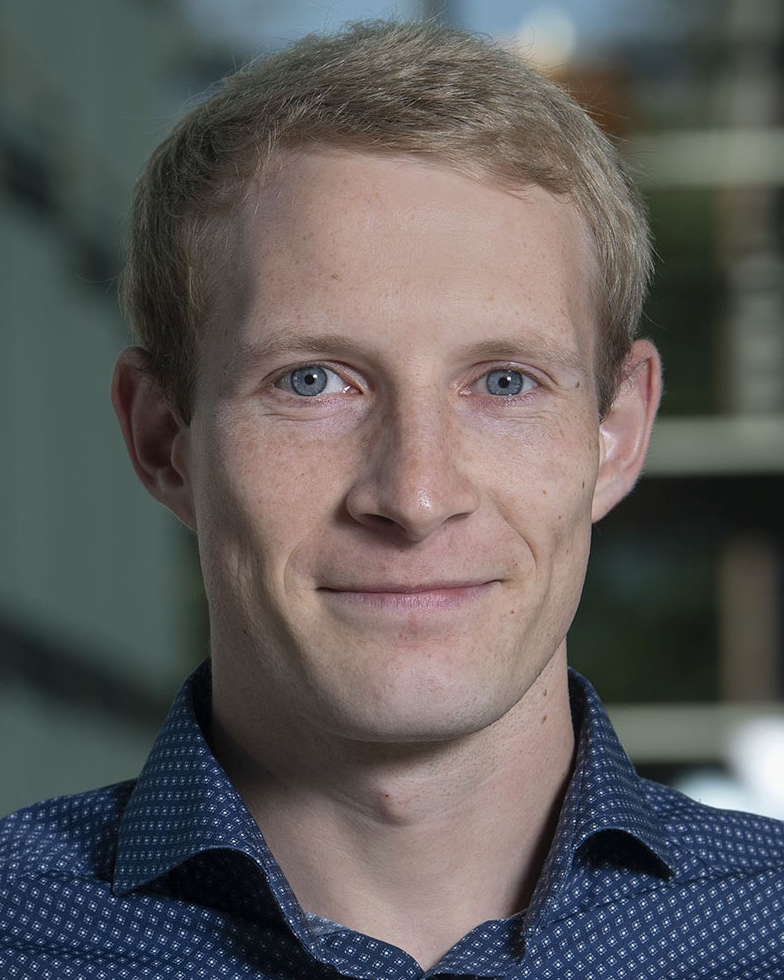}}]{Anders E. Kalør}  (Member, IEEE)  received the B.Sc. and M.Sc. degrees in computer engineering and the Ph.D. degree in wireless communications from Aalborg University, Denmark, in 2015, 2017, and 2022, respectively. He is currently a Project Assistant Professor at the Department of Information and Computer Science, Keio University, Japan. He was a postdoctoral researcher at The University of Hong Kong from 2022 to 2023, and at Aalborg University from 2023 to 2024, supported by an International Postdoc grant from the Independent Research Fund Denmark (IRFD). He was awarded the Spar Nord Foundation Research Award for his Ph.D. project in 2023. His current research interests include communication theory and the intersection between wireless communications, machine learning, and information theory for IoT.
\end{IEEEbiography}

\begin{IEEEbiography}[{\includegraphics[width=1in,height=1.25in, clip,keepaspectratio]{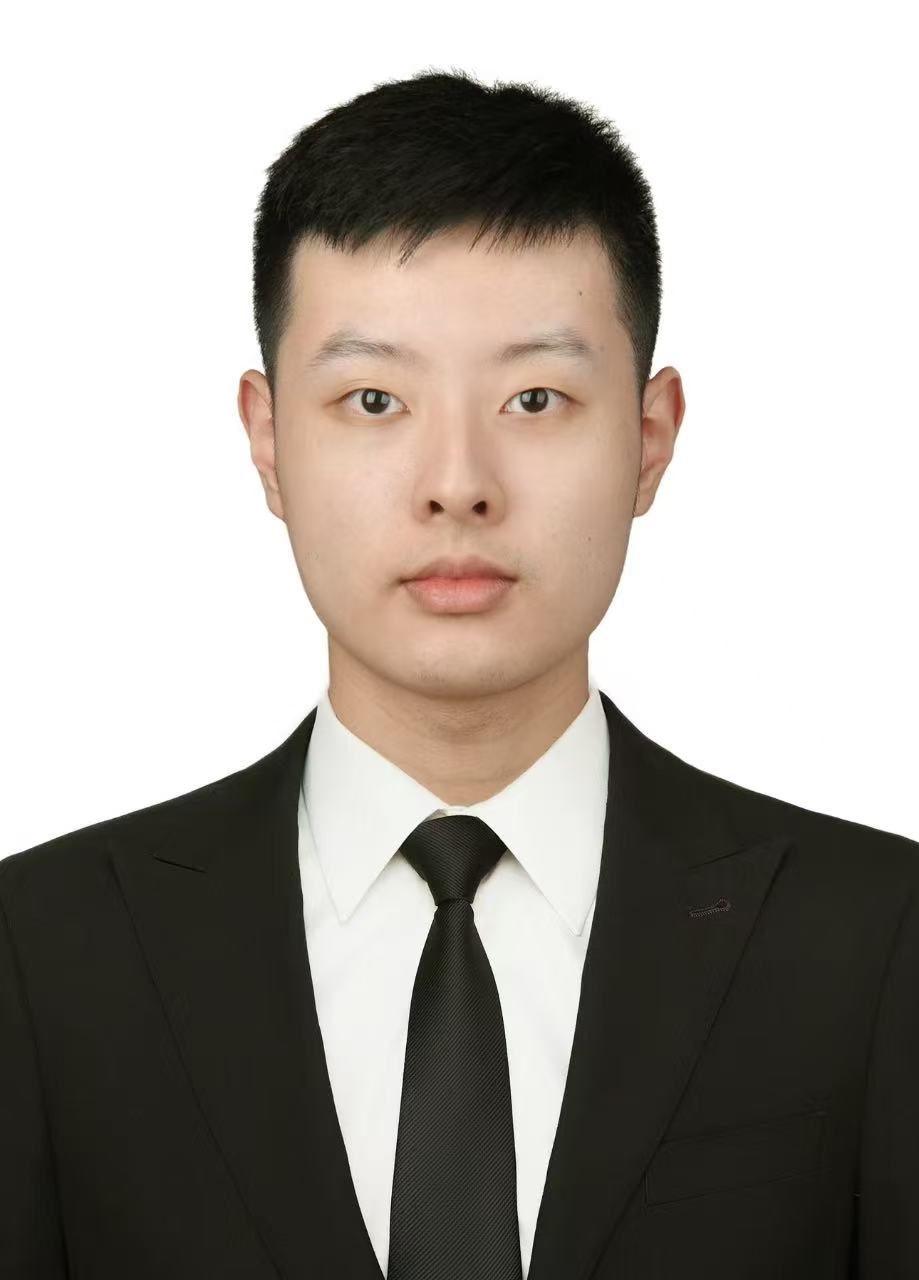}}]{You Zhou} (Graduate Student Member, IEEE) received the B.Eng. in Electrical Engineering from the University of Wisconsin Madison (UW), USA in 2021. He is currently pursuing the Ph.D. degree with the Department of Electrical and Electronic Engineering, The University of Hong Kong, Hong Kong. He is a recipient of Hong Kong PhD Fellowship (HKPF). His research interests include edge inference and 6G wireless communication system design.
\end{IEEEbiography}

\begin{IEEEbiography}[{\includegraphics[width=1in,height=1.25in, clip,keepaspectratio]{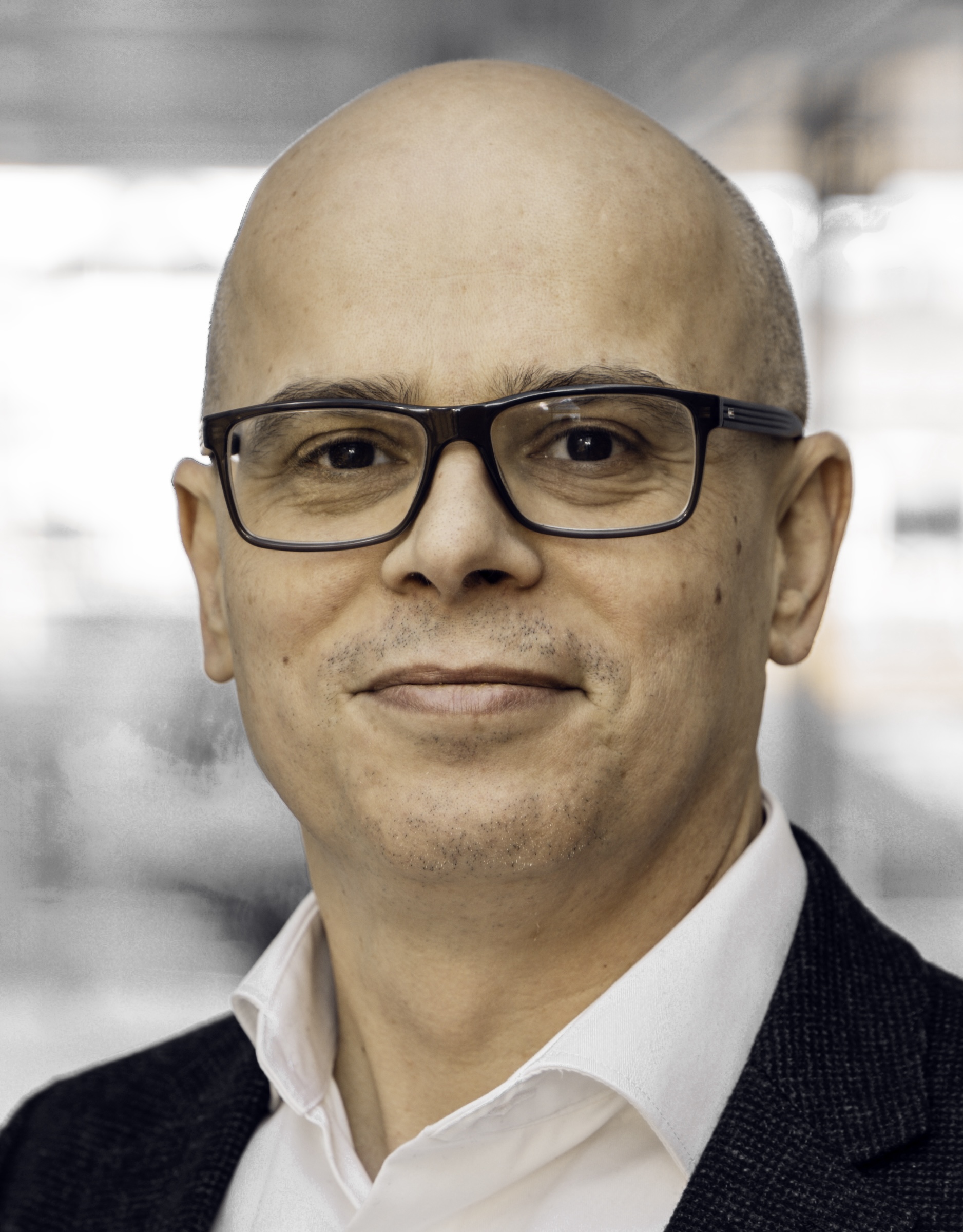}}]{Petar Popovski} (Fellow, IEEE)  is a Professor at Aalborg University, where he is the director of the Center of Excellence CLASSIQUE (Classical Communication in the Quantum Era). He holds a position of a Visiting Excellence Chair at the University of Bremen and a Visiting Professor at the University of Sts. Cyril and Methodius (UKIM) in Skopje. He received his Dipl.-Ing (1997) and M. Sc. (2000) degrees in communication engineering from UKIM and the Ph.D. degree (2005) from Aalborg University. He is a Fellow of the IEEE. He received technical recognition/achievement award from multiple Technical Committees of the IEEE Communications Society: Smart Grid Communications (2019), Wireless Communications (2024), and Satellite and Space Communications (2025). In addition, he received ERC Consolidator Grant (2015), the Danish Elite Researcher award (2016), IEEE Fred W. Ellersick prize (2016), IEEE Stephen O. Rice prize (2018), the Danish Telecommunication Prize (2020), and Villum Investigator Grant (2021). He was a Member at Large at the Board of Governors in IEEE Communication Society and Chair of the IEEE Communication Theory Technical Committee. His research interests are in the area of wireless communication and communication theory. He authored the book ``Wireless Connectivity: An Intuitive and Fundamental Guide'', published by Wiley in 2020. He is currently an Editor-in-Chief of IEEE JOURNAL ON SELECTED AREAS IN COMMUNICATIONS.
\end{IEEEbiography}

\begin{IEEEbiography}[{\includegraphics[width=1in,height=1.25in, clip,keepaspectratio]{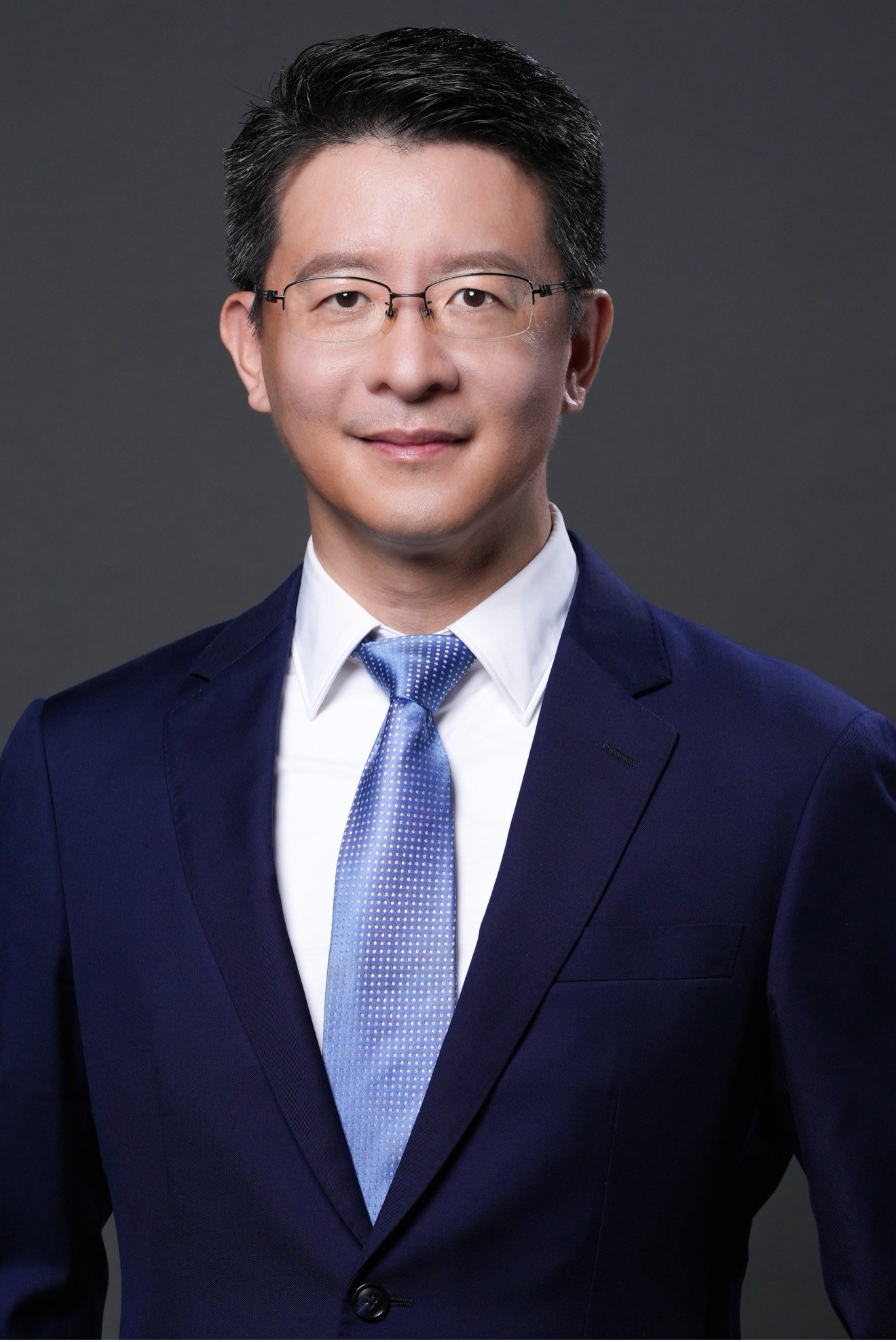}}]{Kaibin Huang} (Fellow, IEEE) received the B.Eng. and M.Eng. degrees from the National University of Singapore and the Ph.D. degree from The University of Texas at Austin, all in electrical engineering. He is the Philip K H Wong Wilson K L Wong Professor in Electrical Engineering and the Department Head at the Dept. of Electrical and Electronic Engineering, The University of Hong Kong (HKU), Hong Kong. His work was recognized with seven Best Paper awards from the IEEE Communication Society. He is a member of the Engineering Panel of Hong Kong Research Grants Council (RGC) and a RGC Research Fellow (2021 Class). He has served on the editorial boards of five major journals in the area of wireless communications and co-edited ten journal special issues. He has been active in organizing international conferences such as the 2014, 2017, and 2023 editions of IEEE Globecom, a flagship conference in communication. He has been named as a Highly Cited Researcher by Clarivate in the last six years (2019-2024) and an AI 2000 Most Influential Scholar (Top 30 in Internet of Things) in 2023-2024. He was an IEEE Distinguished Lecturer. He is a Fellow of the IEEE and the U.S. National Academy of Inventors. 

\end{IEEEbiography}

\end{document}